\documentclass[a4paper,10pt]{amsart}

\usepackage[colorlinks=printing]{hyperref}
\usepackage{amssymb,amsmath,esint,graphicx,subfigure}

\newtheorem{theorem}{Theorem}[section]
\newtheorem{lemma}[theorem]{Lemma}
\newtheorem{proposition}[theorem]{Proposition}
\newtheorem{corollary}[theorem]{Corollary}
\newtheorem{definition}{Definition}[section]

\theoremstyle{remark} \newtheorem{remark}[theorem]{Remark} \theoremstyle{definition} \newtheorem{example}{Example}[section]

\numberwithin{equation}{section}

\newcommand{\Ref}[1]{(\ref{#1})}

\begin{document}

\title[DG method for fractal conservation laws]{The discontinuous Galerkin method\\ for fractal conservation laws}

\author[{S.~Cifani}]{{Simone Cifani}}
\address[Simone Cifani]{\\ Department of Mathematics\\ Norwegian University of Science and Technology (NTNU)\\
 N-7491 Trondheim, Norway}
 \email[]{simone.cifani\@@math.ntnu.no}
 \urladdr{http://www.math.ntnu.no/\~{}cifani/}

\author[{E.R.~Jakobsen}]{{Espen R. Jakobsen}}
\address[Espen R. Jakobsen]{\\ Department of Mathematics\\ Norwegian University of Science and Technology (NTNU)\\
 N-7491 Trondheim, Norway}
\email[]{erj\@@math.ntnu.no}
\urladdr{http://www.math.ntnu.no/\~{}erj/}

\author[{K.H.~Karlsen}]{{Kenneth H. Karlsen}}
\address[Kenneth H. Karlsen]{\\ Centre of Mathematics for
Applications (CMA)\\ Department of Mathematics\\ University of Oslo\\ P.O. Box 1053, Blindern\\ N-0316 Oslo, Norway} \email[]{kennethk\@@math.uio.no}
\urladdr{http://www.math.uio.no/\~{}kennethk/}


\thanks{This research was supported by the Research Council of Norway (NFR) through the project "Integro-PDEs:
Numerical methods, Analysis, and Applications to Finance". The work of K.~H.~Karlsen was also supported through a NFR Outstanding Young Investigator Award.
This article was written as part of the international research program on Nonlinear Partial Differential Equations at the Centre for Advanced Study at the
Norwegian Academy of Science and Letters in Oslo during the academic year 2008--09.}

\begin{abstract}
  We propose, analyze, and demonstrate a discontinuous Galerkin method
  for fractal conservation laws. Various stability estimates are
  established along with error estimates for regular solutions of
  linear equations. Moreover, in the nonlinear case and whenever piecewise
  constant elements are utilized, we prove a rate of convergence
  toward the unique entropy solution. We present numerical results for
  different types of solutions of linear and nonlinear fractal
  conservation laws.
\end{abstract}

\maketitle

\section{Introduction}
We consider the fractional (also called \emph{fractal}) conservation law
\begin{equation}\label{1}
\left\{
\begin{array}{ll}
\partial_{t}u(x,t)+\partial_{x}f(u(x,t))=g_{\lambda}[u(x,t)]&(x,t)\in Q_T:=\mathbb{R}\times(0,T),\\
u(x,0)=u_{0}(x)&x\in\mathbb{R},
\end{array}
\right.
\end{equation}
where $f$ is a Lipschitz continuous function and $g_{\lambda}$ is the nonlocal fractional Laplace operator
$-(-\partial_{x}^{2})^{\lambda/2}$ for some $\lambda\in(0,1)$. This operator can be formally defined by Fourier transform as
\begin{align}\label{fourier}
\widehat{g_{\lambda}[\varphi(x)]}(\xi)=-|\xi|^{\lambda}\hat{\varphi}(\xi)
\end{align}
or, equivalently, by a singular integral (cf.~\cite{Droniou/Imbert,Landkof}) as
\begin{equation*}
\begin{split}
g_\lambda[\varphi(x)]=c_{\lambda}\int_{|z|>0}\frac{\varphi(x+z)-\varphi(x)}{|z|^{1+\lambda}}dz
\end{split}
\end{equation*}
for some $c_\lambda>0$. For sake of brevity, we often write $g$ instead of $g_\lambda$ in the following.

Nonlocal partial differential equations appear in different areas of engineering and sciences. For example, the linear nonlocal partial differential equation
\begin{equation}\label{lin}
\partial_{t}u-\partial_{x}^{2}u-\partial_{x}u+u=g_{\lambda}[u]
\end{equation}
is a nonlocal generalizations of the famous Black-Scholes' equation in finance \cite{Cont/Tankov}, and has received a lot of attention in the last
decade. In recent years, attention has also been given to nonlinear nonlocal equations like
\begin{equation}\label{fracbur}
\partial_{t}u+u\partial_{x}u=g_{\lambda}[u],
\end{equation}
known as the fractional Burgers' equation. Equation \eqref{fracbur} finds application in certain models of detonation of gases (cf.~\cite{Matalon})
characterized by an anomalous diffusive behavior which can be described by means of the fractional Laplacian. We refer the reader to
\cite{Alibaud,Alibaud/Droniou/Vovelle,Droniou}, and the references therein, for further applications in hydrodynamics, molecular biology, semiconductor growth
and dislocation dynamics.

Many authors, see
\cite{Alibaud,Alibaud/Droniou/Vovelle,Biler/Funaki/Woyczynski,Biler/Karch/Woyczynki,Bossy/Jourdain,Brandolese/Karch,Droniou/Imbert,Karch/Miao/Xu}, have
contributed to settle issues like well-posedness and regularity of solutions for the fractional conservation law \eqref{1}. In the case $\lambda\in(1,2)$,
\eqref{1} is the natural nonlocal generalization of the viscous conservation law $\partial_{t}u+\partial_{x}f(u)=\partial_{x}^{2}u$. Such equations turn a
merely bounded initial datum into a unique stable smooth solution (cf.~\cite{Droniou/Gallouet/Vovelle}). The case $\lambda\in(0,1)$ is more delicate. Alibaud's
entropy formulation is needed to guarantee well-posedness \cite{Alibaud}, and the solutions may develop shocks in finite time \cite{Alibaud/Droniou/Vovelle};
the diffusion is no longer strong enough to counterbalance the convection, and equation \eqref{1} fails to regularize the initial datum. In the critical case
$\lambda=1$, Alibaud's entropy formulation is still needed to ensure well-posedness, however, solutions should be smooth as in the case $\lambda\in(1,2)$ --
see Kiselev \emph{et al.~}\cite{Kiselev/Nazarov/Shterenberg} for the case of the fractional Burgers' equation.

A vast literature is available on numerical methods for nonlocal linear equations like \eqref{lin}. The interested reader could see, for example,
\cite{Almendral/Oosterlee,Briani/LaChioma/Natalini,Briani/Natalini,Briani/Natalini/Russo,Cont/Ekaterina,DHalluin/Forsyth/Vetzal,Matache/Schwab/Wihler}.
However, numerical methods for nonlocal nonlinear equations like \eqref{1} are far from being abundant. Dedner \emph{et al.}~introduced in \cite{Dedner/Rohde}
a general class of differences methods for a nonlinear nonlocal equation similar to \eqref{1} coming from a specific problem in radiative hydrodynamics.
Droniou \cite{Droniou} was the first to analyze a general class of difference methods for \eqref{1}, he proved convergence toward Alibaud's entropy solution, but
produced no results regarding the rate of convergence of his methods.

In this paper we study a discontinuous Galerkin (DG) approximation of \Ref{1}. The DG method is a well established numerical method for the pure
conservation law $\partial_{t}u+\partial_{x}f(u)=0$. Some of the important features of this method are stability and high-order accuracy. Moreover, when
piecewise constant elements are used, the DG method reduces to a conservative monotone difference method (cf.~\cite{Holden/Risebro}) which converges to the
entropy solution with rate $1/2$ (cf.~the well known results of Kuznetsov \cite{Kuznetsov}).
For a detailed presentation of the DG method for pure conservation laws, we refer to Cockburn \cite{Cockburn}.

In this paper we propose a DG approximation of \Ref{1} in the case $\lambda\in(0,1)$, and prove that we retain the main features of the DG method in our
nonlocal setting. We show $L^2$-stability, and prove high-order accuracy for linear equations. Moreover, when piecewise constant elements are used, we derive
two fully discrete numerical methods, an implicit-explicit method as in \cite{Droniou} and a fully explicit one, and prove convergence toward a BV entropy
solution of \eqref{1} (cf.~Definition \ref{entropy} below) with a certain rate. For the implicit-explicit method we prove convergence with rate $1/2$ while for
the fully explicit one we prove convergence with a lower rate, $\min\{1/2,1-\lambda\}$. To prove the rate of convergence, we generalize the Kuznetsov argument
\cite{Kuznetsov} to our nonlocal setting, and, as a byproduct, we obtain the following theoretical result: Alibaud's entropy formulation and the BV
entropy formulation are equivalent whenever the initial datum is integrable and of bounded variation.

Finally, several numerical experiments have been performed to illustrate the developed theory. Among other things, we are able to reproduce the theoretical
results (absence of smoothing effect due to persistence of discontinuities and formations of shocks) obtained in
\cite{Alibaud/Droniou/Vovelle,Kiselev/Nazarov/Shterenberg} for the fractional Burgers' equation.

\section{A semidiscrete DG method}
Let us introduce the space grid $x_{i}=i\Delta x$, $i\in\mathbb{Z}$, and let us label $I_{i}=(x_{i},x_{i+1})$. We call $P^{k}(I_{i})$ the set of polynomials of
degree at most $k\in\{0,1,2,\ldots\}$ with support on the interval $I_{i}$, and consider the Legendre polynomials (cf.~\cite{Cockburn} for details)
$$\{\varphi_{0,i},\varphi_{1,i},\ldots,\varphi_{k,i}\},\ \varphi_{j,i}\in P^{j}(I_{i})\text{ for all $j=0,\ldots,k$.}$$
Each $\varphi\in P^{k}(I_{i})$ is a linear combination of the functions $\{\varphi_{0,i},\varphi_{1,i},\ldots,\varphi_{k,i}\}$.

If we multiply \eqref{1} by an arbitrary $\varphi\in P^k(I_i)$, integrate over the interval $I_{i}$, integrate by parts, and replace the flux $f$ by a
numerical flux $F$, we get
\begin{equation}\label{14}
\begin{split}
\int_{I_{i}}u_{t}\varphi-\int_{I_{i}}f(u)\varphi_{x}+F(u_{i+1})\varphi(x_{i+1}^{-})-F(u_{i})\varphi(x_{i}^{+}) =\int_{I_{i}}g[u]\varphi.
\end{split}
\end{equation}
As usual for DG methods, the numerical flux $F(u_{i})=F(u(x_{i}^{-}),u(x_{i}^{+}))$ satisfies the following assumptions:
\begin{itemize}
\item[\emph{A1}:]$F$ is Lipschitz continuous on $\mathbb{R}\times\mathbb{R}$,
\item[\emph{A2}:]$F(a,a)=f(a)$ for all $a\in\mathbb{R}$,
\item[\emph{A3}:]$F$ is non-decreasing with respect to its first variable,
\item[\emph{A4}:]$F$ is non-increasing with respect to its second variable.
\end{itemize}
The goal is to find a function $\tilde{u}:\mathbb{R}\times[0,T]\rightarrow\mathbb{R}$,
\begin{equation}\label{12}
\tilde{u}(x,t)=\sum_{i\in\mathbb{Z}}\sum_{p=0}^{k}U_{p,i}(t)\varphi_{p,i}(x),
\end{equation}
which satisfies \eqref{14} for all $\varphi\in P^{k}(I_i)$, $i\in\mathbb{Z}$. Let us fix $\varphi(x)=\sum_{q=0}^{k}\alpha_{q,i}\varphi_{q,i}(x)$, and plug
\eqref{12} into \eqref{14} to get
\begin{equation*}
\begin{split}
&\sum_{q=0}^{k}\alpha_{q,i}\Bigg(\frac{\Delta x}{2q+1}\frac{d}{dt}U_{q,i}\Bigg)\\
&=\sum_{q=0}^{k}\alpha_{q,i}\Bigg(\int_{I_{i}}f(\tilde{u})\frac{d}{dx}\varphi_{q,i}
+(-1)^{q}F(\tilde{u}_{i})-F(\tilde{u}_{i+1})+\int_{I_{i}}g[\tilde{u}]\varphi_{q,i}\Bigg)
\end{split}
\end{equation*}
where $F(\tilde{u}_{i})=F(\sum_{p=0}^{k} U_{p,i-1},\sum_{p=0}^{k} U_{p,i}(-1)^{p})$. To derive the above expression we have used some well known properties of
the Legendre polynomials: for all $i\in\mathbb{Z}$,
\begin{equation*}
\begin{split}
\int_{I_{i}}\varphi_{p,i}\varphi_{q,i}dx&=\left\{
\begin{array}{cl}
\frac{\Delta x}{2q+1}&\text{for }p=q\\
0&\text{otherwise}
\end{array}
\right.,\ \varphi_{p,i}(x_{i+1}^{-})=1\text{ and }\varphi_{p,i}(x_{i}^{+})=(-1)^{p},
\end{split}
\end{equation*}
where we have denoted with $\varphi(x_i^{+}),\varphi(x_i^{-})$ the (right and left) limits of $\varphi(s)$ as $s\rightarrow x_i$.
The semidiscrete method (i.e., discrete in space and continuous in time) we study is the following: for all $q=0,\ldots,k$ and $i\in\mathbb{Z}$,
\begin{equation}\label{semidicrete_method}
\left\{\begin{array}{rl} \frac{\Delta x}{2q+1}\frac{d}{dt}U_{q,i}\!\!\!&=\int_{I_{i}}f(\tilde{u})\frac{d}{dx}\varphi_{q,i}
+(-1)^{q}F(\tilde{u}_{i})-F(\tilde{u}_{i+1})+\int_{I_{i}}g[\tilde{u}]\varphi_{q,i},\\[0.3cm]
U_{q,i}(0)\!\!\!&=\frac{2q+1}{\Delta x}\int_{I_{i}}u_{0}(x)\varphi_{q,i}(x)dx.
\end{array}\right.
\end{equation}

\section{Nonlinear $L^2$-stability and convergence in the linear case}
Let $V^{k}:=\{u:u|_{I_{i}}\in P^{k}(I_{i})\text{ for all }i\in\mathbb{Z}\}$ be the space of piecewise polynomials, and let $H^{\lambda/2}(\mathbb{R})$
be the fractional Sobolev space with norm
\begin{equation*}
\begin{split}
\|u\|_{H^{\lambda/2}(\mathbb{R})}^{2}:=\|u\|_{L^{2}(\mathbb{R})}^{2}+|u|_{H^{\lambda/2}(\mathbb{R})}^{2}\text{ and }
|u|_{H^{\lambda/2}(\mathbb{R})}^{2}:=\int_{\mathbb{R}}\int_{\mathbb{R}}\frac{[u(z)-u(x)]^{2}}{|z-x|^{1+\lambda}}dzdx.
\end{split}
\end{equation*}
Let us note that the space $H^{\lambda/2}(\mathbb{R})$ also contains {\em discontinuous}
functions (cf.~\cite[Lemma $6.5$]{Folland}). Moreover, let us denote with $H^{-\lambda/2}(\mathbb{R})$
the dual space of $H^{\lambda/2}(\mathbb{R})$, and let us point out that, as shown in the proof of Corollary \ref{023} below, $g[u]\in
H^{-\lambda/2}(\mathbb{R})$ whenever $u\in H^{\lambda/2}(\mathbb{R})$. In the following, all the integrals of the form $\int_{\mathbb{R}}g[u]v$,
where the functions $u,v\in H^{\lambda/2}(\mathbb{R})$, should be interpreted as the pairing $\langle g[u],v \rangle$ between $H^{\lambda/2}(\mathbb{R})$ and its dual.

\begin{theorem}\label{24}\emph{(Stability)} If also $f(0)=0$, then any
  solution $\tilde{u}$ of \eqref{semidicrete_method} belonging to
  $C^1([0,T];H^{\lambda/2}(\mathbb R))$ is $L^2$-stable:
\begin{equation*}
\begin{split}
\|\tilde{u}(\cdot,T)\|_{L^{2}(\mathbb{R})}^2+c_\lambda\int_0^T|\tilde{u}(\cdot,t)|_{H^{\lambda/2}(\mathbb{R})}^{2}dt\leq\|u_{0}\|_{L^{2}(\mathbb{R})}^2.
\end{split}
\end{equation*}
\end{theorem}
The above result generalizes a well known result for the DG method for pure conservation laws (cf.~\cite[Proposition 2.1 and Theorem 4.2]{Cockburn} for
details).
\begin{proof}
  By construction, $\tilde{u}(\cdot,t)$ satisfies \eqref{14} for all
  test functions $\varphi\in P^k(I_i)$. Let us choose the test function
  $\varphi=\tilde{u}(\cdot,t)$, sum over $i\in\mathbb{Z}$, rearrange
  the terms in the sum and integrate over time to get
\begin{equation*}
\begin{split}
\int_{Q_T}\tilde{u}_{t}\tilde{u}=
\int_0^T\sum_{i\in\mathbb{Z}}\Big[F(\tilde{u}_{i})(\tilde{u}(x_{i}^{+})-\tilde{u}(x_{i}^{-}))+\int_{I_{i}}f(\tilde{u})\tilde{u}_{x}\Big]
+\int_{Q_T}g[\tilde{u}]\tilde{u}.
\end{split}
\end{equation*}
Due to the assumptions made (including \emph{A1-A4}), each term in the above expression is well defined. The first term is clear while the last term is well
defined by Corollary \ref{023}. The remaining terms makes sense for all functions in $V^k\cap L^2(\mathbb{R})$ since the point values are well defined. To see
this, note that, since $f(0)=0$ and $f$ is Lipschitz continuous, $f(\tilde u)$ and $F(\tilde u)$ belongs to $L^2(\mathbb{R})$ since $\tilde u$ does. We can
then conclude, using the Cauchy-Schwarz inequality, if the function $v$, $v=\tilde u_x$ in $\cup_{i\in\mathbb{Z}} I_i$,
belongs to $L^2(\mathbb{R})$ (this is the regular part
of the distribution $\tilde u_x$). But this again is an easy consequence of the regularity of the Legendre polynomials $\varphi_{p,i}$
and their othogonality which implies that
$$\sum_{i\in\mathbb{Z}}\sum_{p=0}^{k} c_{p,i}(t)\varphi_{p,i}(x)\in L^2(Q_T)\ \text{if and
  only if}\ \sum_{i\in\mathbb{Z}}\sum_{p=0}^{k}\int_0^T c_{p,i}^2(t)dt <\infty.$$
Let us now prove stability. Since
\begin{equation*}
\begin{split}
\int_{I_{i}}f(u)u_{x}=\int_{I_{i}}(\int^{u(x)}f)_{x}=\int^{u(x_{i+1}^{-})}f-\int^{u(x_{i}^{+})}f,
\end{split}
\end{equation*}
we find that
\begin{equation*}
\begin{split}
\int_{Q_T}\tilde{u}_{t}\tilde{u}
=\int_0^T\sum_{i\in\mathbb{Z}}\Big[F(\tilde{u}_{i})(\tilde{u}(x_{i}^{+})-\tilde{u}(x_{i}^{-}))-\int_{\tilde{u}(x_{i}^{-})}^{\tilde{u}(x_{i}^{+})}f(x)dx\Big]
+\int_{Q_T}g[\tilde{u}]\tilde{u}.
\end{split}
\end{equation*}
It is well known that a flux satisfying \emph{A2}-\emph{A4} is an E-flux (cf.~\cite{Cockburn}), i.e.
\begin{equation*}
F(\tilde{u}_{i})(\tilde{u}(x_{i}^{+})-\tilde{u}(x_{i}^{-}))-\int_{\tilde{u}(x_{i}^{-})}^{\tilde{u}(x_{i}^{+})}f(x)dx\leq0\text{ for all $i\in\mathbb{Z}$.}
\end{equation*}
Thus, by Corollary \ref{023},
\begin{equation*}
\begin{split}
\frac{1}{2}\|\tilde{u}(\cdot,T)\|_{L^{2}(\mathbb{R})}^{2}
+\frac{c_\lambda}{2}\int_0^T\int_{\mathbb{R}}\int_{\mathbb{R}}\frac{(\tilde{u}(z,t)-\tilde{u}(x,t))^{2}}{|z-x|^{1+\lambda}}dzdxdt
\leq\frac{1}{2}\|\tilde{u}_{0}\|_{L^{2}(\mathbb{R})},
\end{split}
\end{equation*}
and the proof is complete.
\end{proof}

In the linear case, equation \eqref{1} reduces to
\begin{equation}\label{5}
\partial_{t}u+c\partial_{x}u=g[u]
\end{equation}
where $c\in\mathbb{R}$. Let us prove the following result.
\begin{proposition}
Let $u_{0}\in H^{k+1}(\mathbb{R})$, $k\geq0$. Then, there exists a unique function $u\in H^{k+1}(Q_T)$ which solves \eqref{5}. Moreover,
\begin{equation}\label{stability}
\|u(\cdot,t)\|_{H^{k+1}(\mathbb{R})}\leq\|u_{0}\|_{H^{k+1}(\mathbb{R})}.
\end{equation}
\end{proposition}
\begin{proof}
Since \eqref{5} is linear, its Fourier transform, $\partial_{t}\hat{u}+i\xi c\hat{u}=-|\xi|^{\lambda}\hat{u}$, has solution
\begin{equation*}
\hat{u}(\xi,t)=\hat{u}_{0}(\xi)e^{-(i\xi c+|\xi|^{\lambda})t}.
\end{equation*}
This implies existence plus, using Plancherel theorem, $L^{2}$-stability and uniqueness. $L^{2}$-stability for (weak) higher derivatives can be obtained as
follows: take the derivative of \eqref{5}, repeat the above procedure, and iterate until the $k$-th derivative. Regularity in time can be shown by using
equation \eqref{5} and regularity in space.
\end{proof}

As pointed out by Cockburn \cite{Cockburn}, in the linear case all relevant numerical fluxes (Godunov, Engquist-Osher, Lax-Friedrichs, \emph{etc.}) reduce to
\begin{equation}\label{FLUX}
F(a,b)=\frac{c}{2}(a+b)-\frac{|c|}{2}(b-a).
\end{equation}
We use this flux to prove the following result: the order of the semidiscrete method \eqref{semidicrete_method} increases along with the degree $k$ of the
polynomial basis used.

\begin{theorem}\label{25}\emph{(Convergence)}
Let $u\in H^{k+1}(Q_T)$, $k\geq0$, be a solution of \eqref{5} and $\tilde{u}\in C^1([0,T];H^{\lambda/2}(\mathbb R))$ be a solution
of the semidiscrete method \eqref{semidicrete_method}. Then, there exists a constant $c_{k,T}>0$ such that
\begin{equation*}
\|u(\cdot,T)-\tilde{u}(\cdot,T)\|_{L^{2}(\mathbb{R})}\leq c_{k,T}\Delta x^{k+\frac{1}{2}}.
\end{equation*}
\end{theorem}

The above result, called high-order accuracy, generalizes a well known feature of the DG method for pure conservation laws (cf.~\cite[Theorem
$2.1$]{Cockburn}). We are able to prove this result since, as shown in the proof below, the error due to the local terms ($c_{k,T}\Delta x^{k+1/2}$) is bigger than the
one due to the nonlocal term ($c_{k,T}\Delta x^{k+1-\lambda/2}$).
\begin{proof}
  By construction, for all test functions $\varphi\in V^k\cap
  L^{2}(\mathbb{R})$,
\begin{equation*}
\int_{\mathbb{R}}\tilde{u}_{t}\varphi+\sum_{i\in\mathbb{Z}}\Big[F(\tilde{u}_{i})(\varphi(x_{i}^{-})-\varphi(x_{i}^{+}))-\int_{I_{i}}c\tilde{u}\varphi_{x}\Big]
=\int_{\mathbb{R}}g[\tilde{u}]\varphi.
\end{equation*}
Note that $u$ satisfies the analogous expression
\begin{equation}\label{20}
\int_{\mathbb{R}}u_{t}\varphi+\sum_{i\in\mathbb{Z}}\Big[F(u_{i})(\varphi(x_{i}^{-})-\varphi(x_{i}^{+}))-\int_{I_{i}}cu\varphi_{x}\Big]
=\int_{\mathbb{R}}g[u]\varphi.
\end{equation}
To prove the above relation, let us multiply \eqref{5} by a test function $\varphi$ and integrate over $I_{i}$. Note that, thanks to the $H^{k+1}$-regularity of
$u$, $u$ is continuous (by Sobolev embedding). Thus, since $F$ satisfies assumption \emph{A2}, we get that
\begin{align*}
&\int_{I_{i}}u_{t}\varphi+cu_{x}\varphi-g[u]\varphi\\&=\int_{I_{i}}u_{t}\varphi-\int_{I_{i}}cu\varphi_{x}+F(u_{i+1})\varphi(x_{i+1}^{-})
-F(u_{i})\varphi(x_{i}^{+})-\int_{I_{i}}g[u]\varphi.
\end{align*}
We obtain \eqref{20} by summing over all $i\in\mathbb{Z}$ and rearranging the terms in the sum. Let us introduce the bilinear form
\begin{equation*}
B(e,\varphi):=\int_{\mathbb{R}}e_{t}\varphi+\sum_{i\in\mathbb{Z}}\Big[F(e_{i})(\varphi(x_{i}^{-})-\varphi(x_{i}^{+}))-\int_{I_{i}}ce\varphi_{x}\Big]
-\int_{\mathbb{R}}g[e]\varphi,
\end{equation*}
where $e:=u-\tilde{u}\in H^{\lambda/2}(\mathbb{R})$.
Let us call $\mathbf{u}$ the $L^{2}$-projection of $u$ into $V^k$: i.e., $$\int_{I_{i}}(\mathbf{u}(x)-u(x))\varphi_{ji}(x)dx=0\text{
for all $j=0,\ldots,k$ and $i\in\mathbb{Z}$.}$$ Note that, by Lemma \ref{lem:new}, $\mathbf{u}\in V^k\cap L^{2}(\mathbb{R})$ implies $\mathbf{u}\in
H^{\lambda/2}(\mathbb{R})$. Let us call $\mathbf{e}:=\mathbf{u}-\tilde{u}\in
H^{\lambda/2}(\mathbb{R})$. Since $B(e,\mathbf{e})=0$,
$B(\mathbf{e},\mathbf{e})=B(\mathbf{e}-e,\mathbf{e})=B(\mathbf{u}-u,\mathbf{e})$ or
\begin{align*}
\int_{0}^{T}\int_{\mathbb{R}}\mathbf{e}_{t}\mathbf{e}=&\int_{0}^{T}\int_{\mathbb{R}}(\mathbf{u}-u)_{t}\mathbf{e}
-\int_{0}^{T}\sum_{i\in\mathbb{Z}}\Big[F(\mathbf{e}_{i})(\mathbf{e}(x_{i}^{-})-\mathbf{e}(x_{i}^{+}))-\int_{I_{i}}c\mathbf{e}\mathbf{e}_{x}]\\
&+\int_{0}^{T}\sum_{i\in\mathbb{Z}}\Big[F((\mathbf{u}-u)_{i})(\mathbf{e}(x_{i}^{-})-\mathbf{e}(x_{i}^{+}))-\int_{I_{i}}c(\mathbf{u}-u)\mathbf{e}_{x}\Big]\\
&+\int_{0}^{T}\int_{\mathbb{R}}g[\mathbf{e}]\mathbf{e}-\int_{0}^{T}\int_{\mathbb{R}}g[\mathbf{e}-e]\mathbf{e}.
\end{align*}
Note that, since both $e,\mathbf{e}\in H^{\lambda/2}(\mathbb{R})$, each term in the above expression is well defined
(cf.~the discussion in the proof of Theorem \ref{25}). One can
argue as in \cite[Theorem $2.1$]{Cockburn} to bound the local terms by $c_{k,T}\Delta x^{2k+1}$. Hence,
\begin{align*}
\int_{0}^{T}\int_{\mathbb{R}}\mathbf{e}_{t}\mathbf{e}\leq c_{k,T}\Delta x^{2k+1}
+\int_{0}^{T}\int_{\mathbb{R}}g[\mathbf{e}]\mathbf{e}-\int_{0}^{T}\int_{\mathbb{R}}g[\mathbf{e}-e]\mathbf{e}.
\end{align*}
Let us denote by $\mathcal{I}$ what it is left to estimate on the right-hand side of the above inequality. By Corollary \ref{023}, the $H^{\lambda/2}$-regularity of
both $e,\mathbf{e}$ implies that
\begin{align*}
\mathcal{I}&=\frac{1}{2}\int_{0}^{T}\int_{\mathbb{R}}g[\mathbf{e}]\mathbf{e}+\frac{1}{2}\int_{0}^{T}\int_{\mathbb{R}}g[e]e
-\frac{1}{2}\int_{0}^{T}\int_{\mathbb{R}}g[\mathbf{e}-e](\mathbf{e}-e)\\
&\leq\int_{0}^{T}\|(u-\mathbf{u})(\cdot,t)\|_{H^{\lambda/2}(\mathbb{R})}^{2}dt,
\end{align*}
and, by Lemma \ref{00012bis}, $$\|(u-\mathbf{u})(\cdot,t)\|_{H^{\lambda/2}(\mathbb{R})}^{2}\leq c_{k}\|u(\cdot,t)\|^2_{H^{k+1}(\mathbb{R})}\Delta
x^{2k+2-\lambda}.$$ Thus, using the $H^{k+1}$-stability of $u$, $\int_{0}^{T}\int_{\mathbb{R}}\mathbf{e}_{t}\mathbf{e}\leq c_{k,T}[\Delta x^{2k+1} +\Delta
x^{2k+2-\lambda}]$, and, since $\mathbf{e}(x,0)=0$ and $\|\mathbf{e}\|=\|(u-\tilde{u})-(u-\mathbf{u})\|\geq\|e\|-\|u-\mathbf{u}\|$,
\begin{align*}
\|e(\cdot,T)\|_{L^{2}(\mathbb{R})}^{2}\leq c_{k,T}\Big[\Delta x^{2k+1}+\Delta x^{2k+2-\lambda}+\Delta x^{2k+2}\Big]\leq c_{k,T}\Delta x^{2k+1}.
\end{align*}
\end{proof}

\begin{remark}
Let us prove that a solution $\tilde u\in C^1([0,t];H^{\lambda/2}(\mathbb{R}))$ of the semidiscrete method \eqref{semidicrete_method}
actually exists up to some time $t>0$. We consider the map
$$\tilde u(\cdot,t)\in V^k\cap L^2(\mathbb{R})\rightarrow\mathcal{F}_{\tilde u}^{q,i}(t):=\text{the right-hand side of
\eqref{semidicrete_method},}$$ and call $\mathcal{F}_{\tilde u}(\cdot,t):=\sum_{i\in\mathbb{Z}}\sum_{q=0}^k\mathcal{F}_{\tilde u}^{q,i}(t)\varphi_{q,i}(\cdot).$
Note that, using Corollary \ref{cor:s} (here the assumption $f(0)=0$ is needed),
\begin{equation}\label{aaad}
\tilde u(\cdot,t)\in V^k\cap L^2(\mathbb{R})\Rightarrow\mathcal{F}_{\tilde u}(\cdot,t)\in V^k\cap L^2(\mathbb{R}),
\end{equation}
and, since both $(f,F)$ are Lipschitz continuous, there exists a constant $c>0$ such that, for all $\tilde u,\tilde v\in V^k\in L^2(\mathbb{R})$,
\begin{equation}\label{aaac}
\|(\mathcal{F}_{\tilde u}-\mathcal{F}_{\tilde v})(\cdot,t)\|_{L^2(\mathbb{R})}\leq c\|(\tilde u-\tilde v)(\cdot,t)\|_{L^2(\mathbb{R})}.
\end{equation}
Therefore, thanks to \eqref{aaad} and \eqref{aaac}, an application of the Cauchy-Lipschitz's theorem yields the existence of a time $t>0$ and a unique solution
$$\tilde u\in C^1([0,t];V^k\cap L^{2}(\mathbb{R}))$$ of the semidiscrete method \eqref{semidicrete_method}.
To conclude, note that $V^k\cap L^{2}(\mathbb{R})\subseteq H^{\lambda/2}(\mathbb{R})$ by Lemma \ref{lem:new}.
\end{remark}


\section{Convergence in the nonlinear case}
We study the nonlinear case by using only piecewise constant elements ($k=0$): $$\{\varphi_{0,i},\varphi_{1,i},\ldots,\varphi_{k,i}\}=\{\varphi_{0,i}\}, \quad
\varphi_{0,i}=\mathbf{1}_{I_{i}},$$ where $\mathbf{1}_{I_{i}}:\mathbb{R}\rightarrow\mathbb{R}$ is the indicator function of the interval
$I_{i}=(x_{i},x_{i+1})$. Starting from the semidiscrete method \eqref{semidicrete_method}, we derive two fully discrete methods: an implicit-explicit method
and a fully explicit one. By adapting Kuznetsov's technique \cite{Kuznetsov} to our nonlocal setting, we prove that both methods converge toward a BV entropy
solution of \eqref{1} with a certain rate (cf. Theorem \ref{026}). In Corollary \ref{well-pos}, we show how this result ensures well-posedness for BV entropy
solutions of \eqref{1}. Note that, in the nonlinear case, even when pure conservation laws are considered, no results concerning the rate of convergence are
available for high-order polynomials ($k>0$).

Let us introduce the time grid $t_{n}=n\Delta t$, where $n=\{0,\ldots,N\}$ and $N\Delta t=T$. We discretize the semidiscrete method \eqref{semidicrete_method}
in time to obtain the implicit-explicit method
\begin{equation}\label{implicit}
\left\{
\begin{split}
U_{i}^{n+1}&=U_{i}^{n}-\Delta tD_-F(U_{i}^{n},U_{i+1}^{n})+\Delta tg\langle U^{n+1}\rangle_{i},\\
U_{i}^{0}&=\frac{1}{\Delta x}\int_{I_{i}}u_{0}(x)dx,
\end{split}
\right.
\end{equation}
and the fully explicit one
\begin{equation}\label{explicit}
\left\{
\begin{split}
U_{i}^{n+1}&=U_{i}^{n}-\Delta tD_-F(U_{i}^{n},U_{i+1}^{n})+\Delta tg\langle U^{n}\rangle_{i},\\
U_{i}^{0}&=\frac{1}{\Delta x}\int_{I_{i}}u_{0}(x)dx.
\end{split}
\right.
\end{equation}
Here we have introduce the shorthand notation $D_-F(U_{i}^{n},U_{i+1}^{n}):=\frac{1}{\Delta x}(F(U_{i}^{n},U_{i+1}^{n})-F(U_{i-1}^{n},U_{i}^{n}))$ and the
nonlocal operator
\begin{equation*}
\begin{split}
g\langle U^{n}\rangle_{i}:=\frac{1}{\Delta x}\int_{I_{i}}g[\bar{U}^{n}]dx=\frac{1}{\Delta x}\sum_{j\in\mathbb{Z}}G^i_{j}U_{j}^n,
\end{split}
\end{equation*}
where $G_j^i:=\int_{I_i}g[\mathbf{1}_{I_j}]dx$ (we denote with $\bar{U}^{n}:\mathbb{R}\rightarrow\mathbb{R}$ the step function generated by the grid values
$\{U^{n}_i\}_{i\in\mathbb{Z}}$ such that $\bar{U}^n(x)=U^n_i$ for all $x\in[x_i,x_{i+1})$).

\begin{proposition}\label{007}
For all $(i,j)\in\mathbb{Z}\times\mathbb{Z}$,
\begin{equation*}
\begin{split}
\sum_{k\in\mathbb{Z}}|G^i_k|<\infty,\quad \sum_{k\in\mathbb{Z}}G^i_k=0,\quad G^i_j=G^j_i,\quad G^{i+1}_{j+1}=G_j^i.
\end{split}
\end{equation*}
Moreover, $G^i_j\geq0$ whenever $i\neq j$, while
\begin{equation*}
G^i_i=-d_{\lambda}\Delta x^{1-\lambda},\text{ where
}d_\lambda:=c_\lambda\left(\int_{|z|<1}\frac{dz}{|z|^{\lambda}}+\int_{|z|>1}\frac{dz}{|z|^{1+\lambda}}\right)>0.
\end{equation*}
\end{proposition}
\begin{proof}
See the appendix.
\end{proof}

Let us introduce the CFL condition
\begin{equation}\label{CFLim}
\begin{split}
(F_{1}+F_{2})\frac{\Delta t}{\Delta x}\leq1
\end{split}
\end{equation}
for the implicit-explicit method \eqref{implicit} (here $F_{1},F_{2}$ are the Lipschitz constants of $F$ with respect to its first and second variable) and the
CFL condition
\begin{equation}\label{CFLex}
\begin{split}
(F_{1}+F_{2})\frac{\Delta t}{\Delta x}+d_{\lambda}\frac{\Delta t}{\Delta x^{\lambda}}\leq1
\end{split}
\end{equation}
for the fully explicit method \eqref{explicit}. In what follows, the relevant CFL condition is always assumed to hold.

Let us introduce the time discretization into \eqref{12} as follows:
\begin{equation}\label{interpolation}
\tilde{u}(x,t)=U_{i}^{n}\text{ for all $(x,t)\in[x_i,x_{i+1})\times[t_n,t_{n+1})$}.
\end{equation}
\begin{theorem}\label{025}
Let $u_{0}\in L^{1}(\mathbb{R})\cap BV(\mathbb{R})$. Then, both the implicit-explicit method \eqref{implicit} and the fully explicit method \eqref{explicit} enjoy
the following properties: for all $t\geq 0$,
\begin{itemize}
\item[\emph{i)}]$\|\tilde{u}(\cdot,t)\|_{L^{\infty}(\mathbb{R})}\leq\|u_{0}\|_{L^{\infty}(\mathbb{R})}$,
\item[\emph{ii)}]$\|\tilde{u}(\cdot,t)\|_{L^{1}(\mathbb{R})}\leq\|u_{0}\|_{L^{1}(\mathbb{R})}$,
\item[\emph{iii)}]$|\tilde{u}(\cdot,t)|_{BV(\mathbb{R})}\leq|u_{0}|_{BV(\mathbb{R})}$.
\end{itemize}
Moreover, there exists a constant $c>0$ (whose value is independent of the discretization parameter $\Delta x$) such that, for all $s,t\geq 0$,
\begin{itemize}
\item[\emph{iv)}]$\|\tilde{u}(\cdot,s)-\tilde{u}(\cdot,t)\|_{L^{1}(\mathbb{R})}\leq c(|s-t|+\Delta x).$
\end{itemize}
\end{theorem}
\begin{proof}
We give here the proof for the fully explicit method \eqref{explicit}. The proof for the implicit-explicit method \eqref{implicit} can be found in the
appendix.

Let us point out two consequences of Proposition \ref{007}. In the first place, note that the fully explicit method \eqref{explicit} is conservative. Indeed,
since $\sum_{j\in\mathbb{Z}}|G^i_j|<\infty$ for all $i\in\mathbb{Z}$,
\begin{equation}\label{00017}
\sum_{i\in\mathbb{Z}}\sum_{j\in\mathbb{Z}}|G_{j}^{i}U_{j}^n|=\sum_{j\in\mathbb{Z}}|U_{j}^n|\sum_{i\in\mathbb{Z}}|G_{j}^i|<\infty,
\end{equation}
whenever $\sum_{i\in\mathbb{Z}}|U_{i}^n|<\infty$. Thus, since $\sum_{i\in\mathbb{Z}}G^i_j=0$ for all $j\in\mathbb{Z}$,
\begin{align*}
\sum_{i\in\mathbb{Z}}g\langle{U}^{n}\rangle_{i}=\frac{1}{\Delta x}\sum_{i\in\mathbb{Z}}\sum_{j\in\mathbb{Z}}G^i_{j}U_{j}^n=\frac{1}{\Delta
x}\sum_{j\in\mathbb{Z}}U^n_{j}\sum_{i\in\mathbb{Z}}G^i_{j}=0
\end{align*}
which implies $\sum_{i\in\mathbb{Z}}U_{i}^{n+1}=\sum_{i\in\mathbb{Z}}U_{i}^{n}$. In the second place, note that the fully explicit method \eqref{explicit} is
monotone in view of the CFL condition \eqref{CFLex}.

We are now ready to prove the theorem. Indeed, monotonicity and Proposition \ref{007} ($\sum_{k\in\mathbb{Z}}G^i_k=0$) imply item \emph{i}. The proofs of
items \emph{ii} and \emph{iii} follow, word by word, the ones in \cite[Theorem 3.6]{Holden/Risebro}. Finally, note that, since the numerical flux $F$ is
Lipschitz continuous in both variables, there exists a constant $c>0$ such that
\begin{equation}\label{lip}
\begin{split}
  U_{i}^{n+1}-U_{i}^{n}&=\Delta tD_-F(U_{i}^{n},U_{i+1}^{n})+\frac{\Delta t}{\Delta x}\sum_{j\in\mathbb{Z}}G^i_jU^{n}_j\\
  &\leq c\frac{\Delta t}{\Delta x}|U^n_{i+1}-U^n_{i}|+c\frac{\Delta
    t}{\Delta x}|U^n_{i}-U^n_{i-1}|+\frac{\Delta t}{\Delta
    x}\Big|\sum_{j\in\mathbb{Z}}G^i_jU^{n}_j\Big|.
\end{split}
\end{equation}
Let us multiply both sides of \eqref{lip} by $\Delta x$, and sum over all $i\in\mathbb{Z}$. Since
\begin{equation*}
\begin{split}
  \sum_{i\in\mathbb{Z}}\Big|\sum_{j\in\mathbb{Z}}G^i_jU^{n}_j\Big|\leq\int_{\mathbb{R}}|g[\bar{U}^n]|dx
  &=c_\lambda C\|\bar{U}^n\|_{L^{1}(\mathbb{R})}^{1-\lambda}|\bar{U}^n|_{BV(\mathbb{R})}^{\lambda}\\
  &\leq c_\lambda
  C\|u_0\|_{L^{1}(\mathbb{R})}^{1-\lambda}|u_0|_{BV(\mathbb{R})}^{\lambda}
\end{split}
\end{equation*}
(cf. Lemma \eqref{00012}), we get $\|\bar{U}^{n+1}-\bar{U}^{n}\|_{L^{1}(\mathbb{R})}\leq c\Delta t$ which implies \emph{iv} via \eqref{CFLex}.
\end{proof}

Let us introduce the definition of BV entropy solutions of \eqref{1}. Let $\eta_{k}(u):=|u-k|$, $\eta'_k(u):=\text{sgn}(u-k)$ and
$q_{k}(u):=\eta'_k(u)(f(u)-f(k))$.
\begin{definition}\label{entropy}
A function $u\in L^{\infty}(Q_T)$ is a BV entropy solution of \eqref{1} provided that the following two conditions hold:
\begin{itemize}
\item[\emph{i)}]$u\in C([0,T];L^1(\mathbb{R}))\cap
  L^\infty(0,T;BV(\mathbb{R}))$; 
\item[\emph{ii)}]for all $k\in\mathbb{R}$ and all nonnegative $\varphi\in C_{c}^{\infty}(\overline{Q_T})$,
\begin{equation}\label{BVentropy}
\begin{split}
&\int_{Q_T}\eta_k(u)\varphi_{t}+q_k(u)\varphi_{x}+\eta'_k(u)g[u]\varphi dxdt\\
&+\int_{\mathbb{R}}\eta_{k}(u_{0}(x))\varphi(x,0)dx-\int_{\mathbb{R}}\eta_{k}(u(x,T))\varphi(x,T)dx\geq0.
\end{split}
\end{equation}
\end{itemize}
\end{definition}

The nonlocal term in the above definition is well defined since, by the regularity of $u$, $g[u]$ is integrable over the domain $Q_T$ (this is a consequence of
Lemma \ref{00012}). Note that sufficiently regular solutions of \eqref{1} are solutions according to the above definition while solutions according to the
above definition are weak solutions of \eqref{1} (this can be easily proved by choosing $k$ as the supremum of $|u|$). We refer the reader to Alibaud's paper
\cite{Alibaud} for the precise definition of a weak solution of \eqref{1}.

As already mentioned in the introduction, Alibaud's entropy formulation ensures well-posedness for all bounded initial data. We prove that the BV entropy
formulation is well-posed for all initial data belonging to a smaller set, the set of all integrable functions of bounded variation, and, therefore, Alibaud's
entropy formulation and the BV entropy formulation are equivalent whenever the initial datum lies in this smaller set.

The following lemma generalizes to our nonlocal setting a result due to Kuznetsov \cite{Kuznetsov}, and it is used in the proof of Theorem \ref{026}. Let us
introduce the function $\varphi(x,y,t,s)=\omega_{\epsilon}(x-y)\omega_{\delta}(t-s)$ where $\omega_\alpha\in C_{c}^{\infty}(\mathbb{R})$, $\alpha>0$, can be
built as follows: choose $\omega\in C_{c}^{\infty}(\mathbb{R})$ such that $0\leq\omega\leq1$, $\omega(x)=0$ for all $|x|>1$ and
$\int_{\mathbb{R}}\omega(x)dx=1$; finally, call $\omega_{\alpha}(x):=\omega(x/\alpha)/\alpha$.

\begin{lemma}\label{045}
Let $u_{0}\in L^{1}(\mathbb{R})\cap BV(\mathbb{R})$, $u$ be a BV entropy solution of \eqref{1} and $\tilde{u}:Q_T\rightarrow\mathbb{R}$ be any function such
that items ii-iv in Theorem \ref{025} hold. Let
\begin{equation*}
\begin{split}
\Lambda[u,\varphi,k]:=&\int_{Q_T}\eta_k(u)\varphi_{t}+q_k(u)\varphi_{x}+\eta'_k(u)g[u]\varphi dxdt\\
&+\int_{\mathbb{R}}\eta_{k}(u_{0}(x))\varphi(x,0)dx-\int_{\mathbb{R}}\eta_{k}(u(x,T))\varphi(x,T)dx
\end{split}
\end{equation*}
and $\Lambda_{\epsilon,\delta}[\tilde{u},u]:=\int_{Q_T}\Lambda[\tilde{u},\varphi(\cdot,y,\cdot,s),u(y,s)]dyds$. Then, there exists $c>0$ such that, for all
$\epsilon>0$ and $0<\delta<T$,
$$\|u(\cdot,T)-\tilde{u}(\cdot,T)\|_{L^{1}(\mathbb{R})}\leq c(\epsilon+\delta+\Delta x)-\Lambda_{\epsilon,\delta}[\tilde{u},u].$$
\end{lemma}
\begin{proof}
See the appendix.
\end{proof}

The above Kuznetsov type of lemma allow us to prove the following rates of convergence.

\begin{theorem}\label{026}
Let $u_{0}\in L^{1}(\mathbb{R})\cap BV(\mathbb{R})$ and $u$ be a BV entropy solution of \eqref{1}.
\begin{itemize}
\item[\emph{a)}]If $\tilde{u}$ is the solution of the implicit-explicit method \eqref{implicit},
then there exists a constant $c_T>0$ such that $$\|u(\cdot,T)-\tilde{u}(\cdot,T)\|_{L^{1}(\mathbb{R})}\leq c_{T}\sqrt{\Delta x}.$$
\item[\emph{b)}]If $\tilde{u}$ is the solution of the fully explicit method \eqref{explicit},
then there exists a constant $c_T>0$ such that $$\|u(\cdot,T)-\tilde{u}(\cdot,T)\|_{L^{1}(\mathbb{R})}\leq c_{T}(\sqrt{\Delta x}+\Delta x^{1-\lambda}).$$
\end{itemize}
\end{theorem}

The rate of convergence obtained for the implicit-explicit method \eqref{implicit} generalizes to our nonlocal setting the rate of convergence obtained by
Kuznetsov in \cite{Kuznetsov} for local difference methods for pure conservation laws. We suspect the convergence rate for the fully explicit method
\eqref{implicit} to be suboptimal. Anyway, to the best of our knowledge, no convergence proof for the fully explicit case was available in the literature up to
now (cf.~Droniou \cite{Droniou} for an alternative convergence proof, without convergence rate, for the implicit-explicit case).

\begin{proof}
The plan is to estimate $-\Lambda_{\epsilon,\delta}[\tilde{u},u]$, and, then, use Lemma \ref{045} to conclude.

\emph{Proof for the implicit-explicit method.} Let us introduce the notation $a\wedge b=\min\{a,b\}$, $a\vee b=\max\{a,b\}$, $\eta_{i}^{n}=|U_{i}^{n}-u|$ and
$q_{i}^{n}=f(U_{i}^{n}\vee u)-f(U_{i}^{n}\wedge u)$, where $u=u(y,s)$. Note that $-\Lambda_{\epsilon,\delta}[\tilde{u},u]$ can be rewritten as
\begin{equation}\label{using}
\begin{split}
-\Lambda_{\epsilon,\delta}[\tilde{u},u]
=\int_{Q_T}\Bigg\{&\sum_{n=0}^{N-1}\sum_{i\in\mathbb{Z}}\Big[(\eta_{i}^{n+1}-\eta_{i}^{n})\int_{x_{i}}^{x_{i+1}}\varphi(x,t_{n+1})dx\\
&\qquad\qquad+(q_{i}^{n}-q_{i-1}^{n})\int_{t_{n}}^{t_{n+1}}\varphi(x_{i},t)dt\Big]\\
&-\int_{0}^{T}\int_{\mathbb{R}}\eta'_{u}(\tilde{u})g[\tilde{u}]\varphi dxdt\Bigg\}dyds.
\end{split}
\end{equation}
Indeed, using summation by parts,
\begin{equation*}
\begin{split}
&-\sum_{i\in\mathbb{Z}}\Bigg\{\sum_{n=0}^{N-1}\int_{t_{n}}^{t_{n+1}}\int_{x_{i}}^{x_{i+1}}\eta_{i}^{n}\varphi_{t}(x,t)+q_{i}^{n}\varphi_{x}(x,t)dxdt\\
&\qquad\quad+\eta_{i}^{0}\int_{x_{i}}^{x_{i+1}}\varphi(x,0)dx-\eta_{i}^{N}\int_{x_{i}}^{x_{i+1}}\varphi(x,T)dx\Bigg\}\\
=&-\sum_{i\in\mathbb{Z}}\Bigg\{\sum_{n=0}^{N-1}\eta_{i}^{n}\int_{x_{i}}^{x_{i+1}}\big[\varphi(x,t_{n+1},)-\varphi(x,t_{n},)\big]dx\\
&\qquad\quad+\sum_{n=0}^{N-1}q_{i}^{n}\int_{t_{n}}^{t_{n+1}}\big[\varphi(x_{i+1},t)-\varphi(x_{i},t)\big]dt\\
&\qquad\quad+\eta_{i}^{0}\int_{x_{i}}^{x_{i+1}}\varphi(x,0)dx-\eta_{i}^{N}\int_{x_{i}}^{x_{i+1}}\varphi(x,T)dx\Bigg\}\\
=&\sum_{i\in\mathbb{Z}}\sum_{n=0}^{N-1}\Big[(\eta_{i}^{n+1}-\eta_{i}^{n})\int_{x_{i}}^{x_{i+1}}\varphi(x,t_{n+1})dx
\\&\qquad\qquad+(q_{i}^{n}-q_{i-1}^{n})\int_{t_{n}}^{t_{n+1}}\varphi(x_{i},t)dt\Big].
\end{split}
\end{equation*}
Let us exploit monotonicity to get
\begin{align*}
&U_{i}^{n+1}\vee k\leq U_{i}^{n}\vee k-\Delta tD_-F(U_{i}^{n}\vee k,U_{i+1}^{n}\vee k)+\Delta t\mathbf{1}_{(k,+\infty)}(U_{i}^{n+1})g\langle U^{n+1}\rangle_{i},\\
&U_{i}^{n+1}\wedge k\geq U_{i}^{n}\wedge k-\Delta tD_-F(U_{i}^{n}\wedge k,U_{i+1}^{n}\wedge k)+\Delta t\mathbf{1}_{(-\infty,k)}(U_{i}^{n+1})g\langle
U^{n+1}\rangle_{i}.
\end{align*}
Let us call $Q_{i}^{n}:=F(U_{i}^{n}\vee k,U_{i+1}^{n}\vee k)-F(U_{i}^{n}\wedge k,U_{i+1}^{n}\wedge k)$, and note that, since $|a-b|=a\vee b-a\wedge b$, we can subtract
$U_{i}^{n+1}\wedge k$ from $U_{i}^{n+1}\vee k$ to obtain the cell entropy inequality
\begin{equation}\label{cell_entropy_inequality}
\begin{split}
\eta_{i}^{n+1}-\eta_{i}^{n}+\frac{\Delta t}{\Delta x}(Q_{i}^{n}-Q_{i-1}^{n})-\Delta t\eta'_k(U^{n+1}_i)g\langle U^{n+1}\rangle_{i}\leq0.
\end{split}
\end{equation}
If we plug the above inequality into \eqref{using}, we find that
\begin{equation*}
\begin{split}
-\Lambda_{\epsilon,\delta}[\tilde{u},u]\leq\int_{Q_T}\Bigg\{&\sum_{n=0}^{N-1}\sum_{i\in\mathbb{Z}}
\Big[(q_{i}^{n}-q_{i-1}^{n})\int_{t_{n}}^{t_{n+1}}\varphi(x_{i},t)dt\\
&-\frac{\Delta t}{\Delta x}(Q_{i}^{n}-Q_{i-1}^{n})\int_{x_{i}}^{x_{i+1}}\varphi(x,t_{n+1})dx\Big]\\
&+\Delta t\sum_{n=0}^{N-1}\sum_{i\in\mathbb{Z}}
\eta'_{u}(U^{n+1}_i)g\langle U^{n+1}\rangle_{i}\int_{x_{i}}^{x_{i+1}}\varphi(x,t_{n+1})dx\\
&-\int_{Q_T}\eta'_{u}(\tilde{u})g[\tilde{u}]\varphi dxdt\Bigg\}dyds.
\end{split}
\end{equation*}

Next, the right-hand side of the above inequality needs to be estimated. To this end, let us point out that, as proved in \cite[Example 3.14]{Holden/Risebro},
\begin{equation*}
\begin{split}
\int_{Q_T}\Bigg\{\sum_{n=0}^{N-1}\sum_{i\in\mathbb{Z}}\Big[(q_{i}^{n}-q_{i-1}^{n})\int_{t_{n}}^{t_{n+1}}\varphi(x_{i},t)dt&\\
-\frac{\Delta t}{\Delta x}(Q_{i}^{n}-Q_{i-1}^{n})\int_{x_{i}}^{x_{i+1}}\varphi(x,t_{n+1})dx\Big]\Bigg\}dyds&\leq c_{T}\left(\frac{\Delta
x}{\epsilon}+\frac{\Delta x}{\delta}\right).
\end{split}
\end{equation*}
Let us call $\varphi_i^n:=\frac{1}{\Delta x}\int_{x_{i}}^{x_{i+1}}\varphi(x,t_{n})dx$ and $\bar{\varphi}$ the step function built from $\{\varphi_i^n\}$ by
taking $\bar{\varphi}(x,t)=\varphi_i^n$ for all $(x,t)\in[x_i,x_{i+1})\times[t_n,t_{n+1})$. Moreover, let us call $J$ the term which still needs to be
estimated,
\begin{equation*}
\begin{split}
J:=\int_{Q_T}\Bigg\{&\Delta t\Delta x\sum_{n=0}^{N-1}\sum_{i\in\mathbb{Z}} \eta'_{u}(U^{n+1}_i)g\langle
U^{n+1}\rangle_{i}\bar{\varphi}(x_i,t_{n+1})\\
&-\int_{Q_T}\eta'_{u}(\tilde{u})g[\tilde{u}]\varphi dxdt\Bigg\}dyds.
\end{split}
\end{equation*}
Since $g\langle{U}^{n}\rangle_{i}=\frac{1}{\Delta x}\int_{x_{i}}^{x_{i+1}}g[\bar{U}^{n}]dx$, we can rewrite $J$ as
\begin{equation*}
\begin{split}
J=\int_{Q_T}\Bigg\{\int_{\Delta t}^{T+\Delta t} \int_{\mathbb{R}}\eta'_{u}(\tilde{u})g[\tilde{u}]\bar{\varphi}dxdt
-\int_{Q_T}\eta'_{u}(\tilde{u})g[\tilde{u}]\varphi dxdt\Bigg\}dyds
\end{split}
\end{equation*}
which can be split into $J_1+J_2-J_3$, where
\begin{equation*}
\begin{split}
J_1&:=\int_{Q_T}\Bigg\{
\int_{Q_T}\eta'_{u}(\tilde{u})g[\tilde{u}]\big(\bar{\varphi}-\varphi\big)dxdt\Bigg\}dyds,\\
J_2&:=\int_{Q_T}\Bigg\{\int_{T}^{T+\Delta t}\int_{\mathbb{R}}\eta'_{u}(\tilde{u})g[\tilde{u}]\bar{\varphi}dxdt\Bigg\}dyds,\\
J_3&:=\int_{Q_T}\Bigg\{\int_{0}^{\Delta t}\int_{\mathbb{R}}\eta'_{u}(\tilde{u})g[\tilde{u}]\varphi dxdt\Bigg\}dyds.
\end{split}
\end{equation*}
By Lemma \ref{00012} and Theorem \ref{025}, $g[\tilde{u}]\in L^1(Q_T)$ and, thus, both
\begin{equation*}
\begin{split}
J_2&\leq\int_{T}^{T+\Delta t} \int_{\mathbb{R}}|g[\tilde{u}]|\Bigg\{\int_{Q_T}\bar{\varphi}dyds\Bigg\}dxdt,\\
J_3&\leq\int_{0}^{\Delta t}\int_{\mathbb{R}}|g[\tilde u]|\Bigg\{\int_{Q_T}\varphi dyds\Bigg\}dxdt
\end{split}
\end{equation*}
are of order $\Delta x$ (here, as the the following, we use the CFL condition to pass from $\Delta t$ to $\Delta x$). Moreover,
\begin{equation*}
\begin{split}
J_1\leq\int_{Q_T}|g[\tilde{u}]|\Bigg\{\int_{Q_T}|\bar{\varphi}-\varphi|dyds\Bigg\}dxdt\leq c_{T}\left(\frac{\Delta x}{\epsilon}+\frac{\Delta x}{\delta}\right)
\end{split}
\end{equation*}
since, for all $(x,t)\in Q_T$, there exists a constant $c>0$ such that
\begin{equation}\label{jjj}
\begin{split}
\int_{Q_T}|\bar{\varphi}-\varphi|dyds\leq c\left(\frac{\Delta x}{\epsilon}+\frac{\Delta x}{\delta}\right).
\end{split}
\end{equation}

We now prove \eqref{jjj}. Let us call $\bar{\omega}_\epsilon$ the step function built from $\{\omega_{\epsilon,i}\}$, $\omega_{\epsilon,i}=\frac{1}{\Delta
x}\int_{x_{i}}^{x_{i+1}}\omega_\epsilon(s)ds$, as follows: $\bar{\omega}_\epsilon(x)=\omega_{\epsilon,i}$ for all $x\in[x_i,x_{i+1})$. First, note that
\begin{equation}\label{ijj}
\begin{split}
\int_{\mathbb{R}}|\bar{\omega}_\epsilon(x)-\omega_\epsilon(x)|dx\leq \Delta x|\omega_\epsilon|_{BV(\mathbb{R})}.
\end{split}
\end{equation}
Indeed,
\begin{equation*}
\begin{split}
\int_{\mathbb{R}}|\bar{\omega}_\epsilon(x)-\omega_\epsilon(x)|dx
&=\sum_{i\in\mathbb{Z}}\int_{x_{i}}^{x_{i+1}}|\bar{\omega}_\epsilon(x)-\omega_\epsilon(x)|dx\\
&=\sum_{i\in\mathbb{Z}}\int_{x_{i}}^{x_{i+1}}
\left|\frac{1}{\Delta x}\int_{x_{i}}^{x_{i+1}}\omega_\epsilon(s)ds
-\omega_\epsilon(x)\right|dx\\
&\leq\frac{1}{\Delta x}\sum_{i\in\mathbb{Z}}\int_{x_{i}}^{x_{i+1}}
\int_{x_{i}}^{x_{i+1}}|\omega_\epsilon(s)-\omega_\epsilon(x)|dsdx\\
&\leq\frac{1}{\Delta x}\sum_{i\in\mathbb{Z}}|\omega_\epsilon|_{BV(I_i)}\int_{x_{i}}^{x_{i+1}}
\int_{x_{i}}^{x_{i+1}}dsdx\\
&\leq \Delta x|\omega_\epsilon|_{BV(\mathbb{R})}.
\end{split}
\end{equation*}
Next, we note that for all $(x,t)\in Q_T$,
\begin{equation}\label{kkk}
\begin{split}
\int_{Q_T}|\bar{\varphi}-\varphi|dyds =\int_{0}^{T}\int_{\mathbb{R}}|\bar{\omega}_\epsilon(x-y)\omega_\delta(t_n-s)-\omega_\epsilon(x-y)\omega_\delta(t-s)|dyds,
\end{split}
\end{equation}
where $t_n$ is such that $t\in(t_n,t_{n+1})$. Moreover, using \eqref{ijj},
\begin{equation}\label{se1}
\begin{split}
\int_{\mathbb{R}}|\bar{\omega}_\epsilon(x-y)-\omega_\epsilon(x-y)|dy=\int_{\mathbb{R}}|\bar{\omega}_\epsilon(y)-\omega_\epsilon(y)|dy\leq\Delta
x|\omega_\epsilon|_{BV(\mathbb{R})}=c\frac{\Delta x}{\epsilon}
\end{split}
\end{equation}
while, since $t\in(t_n,t_{n+1})$,
\begin{equation}\label{se}
\begin{split}
\int_{0}^{T}|\omega_\delta(t_n-s)-\omega_{\delta}(t-s)|ds\leq\Delta t|\omega_\delta|_{BV(\mathbb{R})}=c\frac{\Delta x}{\delta}.
\end{split}
\end{equation}
Thanks to the estimates \eqref{se1} and \eqref{se}, an application of the triangular inequality to the right-hand side of \eqref{kkk} yields \eqref{jjj}.

The above estimates ensure that $-\Lambda_{\epsilon,\delta}[\tilde{u},u]\leq c_{T}(\frac{\Delta x}{\epsilon}+\frac{\Delta x}{\delta})$. Therefore, we can use
Lemma \ref{045} to obtain
\begin{equation*}
\begin{split}
\|u(\cdot,T)-\tilde{u}(\cdot,T)\|_{L^{1}(\mathbb{R})}\leq c_{T}\left(\epsilon+\delta+\Delta x+\frac{\Delta x}{\epsilon}+\frac{\Delta x}{\delta}\right).
\end{split}
\end{equation*}
The conclusion follows by setting $\epsilon=\delta=\sqrt{\Delta x}$.

\emph{Proof for the fully explicit method.} Let us exploit monotonicity to get
\begin{align*}
&U_{i}^{n+1}\vee k\leq U_{i}^{n}\vee k-\Delta tD_-F(U_{i}^{n}\vee k,U_{i+1}^{n}\vee k)+\Delta t\mathbf{1}_{(k,+\infty)}(U_{i}^{n+1})g\langle U^{n}\rangle_{i},\\
&U_{i}^{n+1}\wedge k\geq U_{i}^{n}\wedge k-\Delta tD_-F(U_{i}^{n}\wedge k,U_{i+1}^{n}\wedge k)+\Delta t\mathbf{1}_{(-\infty,k)}(U_{i}^{n+1})g\langle
U^{n}\rangle_{i}.
\end{align*}
Proceeding as done in the proof for the implicit-explicit method, we obtain the cell entropy inequality
\begin{equation*}
\begin{split}
\eta_{i}^{n+1}-\eta_{i}^{n}+\frac{\Delta t}{\Delta x}(Q_{i}^{n}-Q_{i-1}^{n})-\Delta t\eta'_k(U^{n+1}_i)g\langle U^{n}\rangle_{i}\leq0.
\end{split}
\end{equation*}
Let us add and subtract $\Delta t\eta'_k(U^{n+1}_i)g\langle U^{n+1}\rangle_{i}$ to the left-hand side of the above inequality, and let us use the fact that the
operator $g\langle\cdot\rangle$ is linear to obtain
\begin{equation*}
\begin{split}
\eta_{i}^{n+1}-\eta_{i}^{n}&+\frac{\Delta t}{\Delta x}(Q_{i}^{n}-Q_{i-1}^{n})\\&-\Delta t\eta'_k(U^{n+1}_i)g\langle U^{n}-U^{n+1}\rangle_{i}-\Delta
t\eta'_k(U^{n+1}_i)g\langle U^{n+1}\rangle_{i}\leq0.
\end{split}
\end{equation*}
If we plug the above inequality into \eqref{using}, we find that
\begin{equation*}
\begin{split}
-\Lambda_{\epsilon,\delta}[\tilde{u},u]\leq\int_{Q_T}\Bigg\{&\sum_{n=0}^{N-1}\sum_{i\in\mathbb{Z}}
\Big[(q_{i}^{n}-q_{i-1}^{n})\int_{t_{n}}^{t_{n+1}}\varphi(x_{i},t)dt\\
&-\frac{\Delta t}{\Delta x}(Q_{i}^{n}-Q_{i-1}^{n})\int_{x_{i}}^{x_{i+1}}\varphi(x,t_{n+1})dx\Big]\\
&+\Delta t\sum_{n=0}^{N-1}\sum_{i\in\mathbb{Z}}
\eta'_{u}(U^{n+1}_i)g\langle U^{n+1}\rangle_{i}\int_{x_{i}}^{x_{i+1}}\varphi(x,t_{n+1})dx\\
&+\Delta t\sum_{n=0}^{N-1}\sum_{i\in\mathbb{Z}}
\eta'_{u}(U^{n+1}_i)g\langle U^{n}-U^{n+1}\rangle_{i}\int_{x_{i}}^{x_{i+1}}\varphi(x,t_{n+1})dx\\
&-\int_{Q_T}\eta'_{u}(\tilde{u})g[\tilde{u}]\varphi dxdt\Bigg\}dyds.
\end{split}
\end{equation*}
The only term left to estimate is
\begin{equation*}
\begin{split}
&\int_{Q_T}\Bigg\{\Delta t\Delta x\sum_{n=0}^{N-1}\sum_{i\in\mathbb{Z}}\eta'_{u}(U^{n+1}_i)g\langle U^{n}
-U^{n+1}\rangle_{i}\bar{\varphi}(x_i,t_{n+1})\Bigg\}dyds\\
&\leq\int_{Q_T}\Bigg\{\Delta t\Delta x\sum_{n=0}^{N-1}\sum_{i\in\mathbb{Z}}|g\langle U^{n} -U^{n+1}\rangle_{i}|\bar{\varphi}(x_i,t_{n+1})\Bigg\}dyds\\
&\leq\int_{Q_T}|g[\tilde{u}(x,t)-\tilde{u}(x,t+\Delta t)]|\Bigg\{\int_{Q_T}\bar{\varphi}dyds\Bigg\}dxdt.
\end{split}
\end{equation*}
Note that, using Lemma \ref{00012} and Theorem \ref{025} (item \emph{iv}), the right-hand side of the above inequality is easily seen to be of order $\Delta
x^{1-\lambda}$.

Finally, using Lemma \ref{045},
\begin{equation*}
\begin{split}
\|u(\cdot,T)-\tilde{u}(\cdot,T)\|_{L^{1}(\mathbb{R})}\leq c_{T}\left(\epsilon+\delta+\Delta x+\Delta x^{1-\lambda}+\frac{\Delta x}{\epsilon}+\frac{\Delta
x}{\delta}\right),
\end{split}
\end{equation*}
and the conclusion follows by setting $\epsilon=\delta=\sqrt{\Delta x}$.
\end{proof}

We conclude this paper by proving the following result which is a consequence of Theorem \ref{026}: the definition of a BV entropy solution of \eqref{1} is
well-posed.
\begin{corollary}\label{well-pos}
Let $u_{0}\in L^{1}(\mathbb{R})\cap BV(\mathbb{R})$. Then, there exists a unique BV entropy solution of \eqref{1}.
\end{corollary}

\begin{proof}
Let us give the proof using the implicit-explicit method \eqref{implicit}. Needless to say, the fully explicit method \eqref{explicit} would also do.

\emph{Uniqueness.} Let us assume that both $u$ and $v$ are BV entropy solutions of \eqref{1}. If we add and subtract the solution of the implicit-explicit
method \eqref{implicit}, we obtain
\begin{equation*}
\begin{split}
\|u(\cdot,T)-v(\cdot,T)\|_{L^{1}(\mathbb{R})}&\leq\|u(\cdot,T)-\tilde{u}(\cdot,T)\|_{L^{1}(\mathbb{R})} +\|v(\cdot,T)-\tilde{u}(\cdot,T)\|_{L^{1}(\mathbb{R})}
\end{split}
\end{equation*}
which, by Theorem \ref{026}, is less than or equal to $c_T\sqrt{\Delta x}+c_T\sqrt{\Delta x}$ for all $\Delta x>0$. Therefore, uniqueness follows.

\emph{Existence.} Using a standard argument (cf., for example, \cite[Theorem 3.8]{Holden/Risebro}), Helly's theorem yields the existence of a subsequence
$\tilde{u}\rightarrow u$ in $L^1_{\mathrm{loc}}(Q_T)$ as $\Delta x\rightarrow0$. Moreover,
$u\in C([0,T];L^1(\mathbb R))\cap L^\infty(0,T;BV(\mathbb{R}))$ by Theorem
\ref{025}. To prove that $u$ satisfies the entropy inequality \eqref{BVentropy}, we start from the cell entropy inequality \eqref{cell_entropy_inequality}. Let us choose a
nonnegative test function $\varphi\in C_c^\infty(\overline{Q_T})$ and call $\varphi_i^n:=\varphi(x_i,t_n)$. If we multiply both sides of \eqref{cell_entropy_inequality}
by $\varphi^n_i\geq0$, sum over $i$ and $n$, and use summations by parts, we find that
\begin{equation*}
\begin{split}
\Delta x\Delta t\sum_{n=1}^{N-1}&\sum_{i\in\mathbb{Z}}\eta_{i}^{n}\frac{\varphi_{i}^{n+1}-\varphi_{i}^{n}}{\Delta t}\\
&+\Delta x\Delta t\sum_{n=0}^{N}\sum_{i\in\mathbb{Z}}\left\{Q_i^n\frac{\varphi_{i+1}^{n}-\varphi_{i}^{n}}{\Delta x}
+\eta'_{k}(U^{n+1}_i)g\langle{U}^{n+1}\rangle_i\varphi_i^{n}\right\}\\
&+\Delta x\sum_{i\in\mathbb{Z}}\left\{\varphi_{i}^{0}\eta_{i}^{0}-\varphi_{i}^{N}\eta_{i}^{N}\right\}\geq0.
\end{split}
\end{equation*}
A standard argument shows that all the local terms in the above expression converge to the ones appearing in the inequality \eqref{BVentropy}, cf.~
e.g.~\cite[Theorem 3.9]{Holden/Risebro}. Let us now consider the nonlocal term. Note that (here $\bar{\varphi}$ is as in the proof of Theorem \ref{026})
\begin{equation*}
\begin{split}
&\Delta x\Delta
t\sum_{n=0}^{N}\sum_{i\in\mathbb{Z}}\eta'_{k}(U^{n+1}_i)g\langle{U}^{n+1}\rangle_i\varphi_i^{n}\\
&=\Delta x\Delta t\sum_{n=0}^{N}\sum_{i\in\mathbb{Z}}\eta'_{k}(U^{n+1}_i)g\langle{U}^{n+1}\rangle_i(\varphi_i^{n}-\varphi_i^{n+1}) +\int_{\Delta t}^{T+\Delta
t}\int_{\mathbb{R}}\eta'_{k}(\tilde{u})g[\tilde{u}]\bar{\varphi}dxdt
\end{split}
\end{equation*}
where, since there exists a constant $c>0$ such that $|\varphi_i^{n}-\varphi_i^{n+1}|\leq c\Delta x$ for all $(i,n)$,
\begin{equation*}
\begin{split}
\left|\Delta x\Delta t\sum_{n=0}^{N}\sum_{i\in\mathbb{Z}}\eta'_{k}(U^{n+1}_i)g\langle{U}^{n+1}\rangle_i(\varphi_i^{n}-\varphi_i^{n+1})\right|\leq c\Delta
x\int_{Q_T}|g[\tilde{u}(x,t+\Delta t)]|dxdt.
\end{split}
\end{equation*}
Since $g[\tilde{u}]\in L^1(Q_T)$, the right-hand side of the above expression is of order $\Delta x$. To conclude, we prove that there exists a subsequence
$\{\tilde u\}$ such that
\begin{equation}\label{hhh}
\begin{split}
\int_{\Delta t}^{T+\Delta t}\int_{\mathbb
  R}\eta'_{k}(\tilde{u})g[\tilde{u}]\bar{\varphi}dxdt\stackrel{\Delta x\rightarrow 0}{\longrightarrow}\int_{Q_T}\eta'_k(u)g[u]\varphi
dxdt
\end{split}
\end{equation}
for a.e.~$k\in\mathbb{R}$. This is a consequence of the dominated convergence theorem since the left hand side integrand converges pointwise a.e.~to the right
hand side integrand. Indeed, first note that $\bar{\varphi}\rightarrow\varphi$  pointwise and that a subsequence $\tilde{u}\rightarrow u$ a.e.~in $Q_T$.
Moreover, for a.e.~$k\in\mathbb{R}$ the measure of $\{(x,t)\in Q_T:u(x,t)=k\}$ is null. This means that $\eta_k'(\tilde{u})\rightarrow\eta'_k(u)$ a.e.~in
$Q_T$, since $\eta'_k$ is continuous on $\mathbb{R}\backslash \{k\}$. Finally, by Theorem \ref{026},
\begin{equation*}
\begin{split}
\int_{Q_T}|g[\tilde{u}-u]|dxdt\leq c\int_{0}^{T}\|\tilde{u}-u\|_{L^1(\mathbb{R})}^{1-\lambda}dt\leq c_T\Delta x^{\frac{1-\lambda}{2}}
\end{split}
\end{equation*}
for all $\Delta x>0$, and hence a subsequence $g[\tilde{u}]\rightarrow g[u]$ a.e.~in $Q_T$. The proof for all $k\in\mathbb{R}$ follows the one given by Droniou
in \cite{Droniou}, and this completes the proof.
\end{proof}


\section{Numerical experiments}
We have implemented the numerical method \eqref{semidicrete_method} in the cases $k=0,1,2$ with fully explicit time discretization. To perform computations, we
have set our numerical solutions to zero outside the region $\Omega=\{(x,t):|x|\leq3/2, t\geq0\}$. In other words, we have computed the value
$U_{p,i}(t_{n+1})$ using only the values $\{U_{p,i}(t_n)\}$, where $x_i\in\Omega$ and $p=0,\ldots,k$. This has been done also at the boundaries $|x|=3/2$.

\begin{remark}\label{rem}
Due to infinite speed of propagation (cf.~\cite{Alibaud}), solutions of \eqref{1} do not have, in general, compact support. Therefore, the use of the region
$\Omega$ introduces an additional error which we have not considered in Theorem \ref{25} and Theorem \ref{026}.
\end{remark}

\begin{figure}[t]
\subfigure[$T=0.5$]{
\includegraphics[scale=0.31]{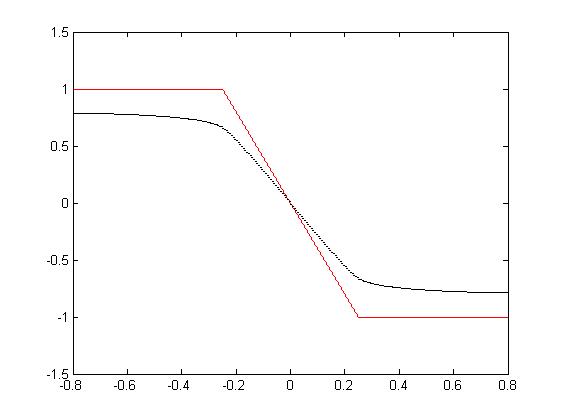}}
\subfigure[$T=1.3$]{
\includegraphics[scale=0.31]{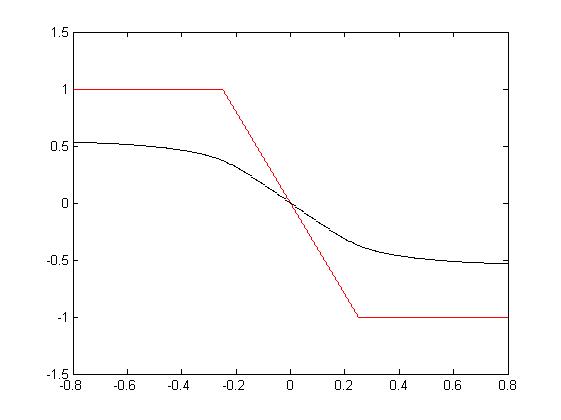}}
\subfigure[$T=0.5$]{
\includegraphics[scale=0.31]{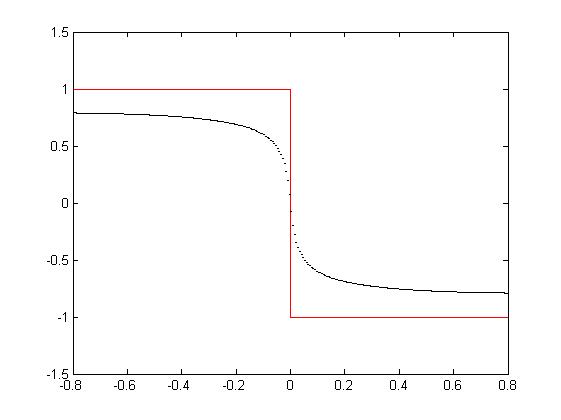}}
\subfigure[$T=1.3$]{
\includegraphics[scale=0.31]{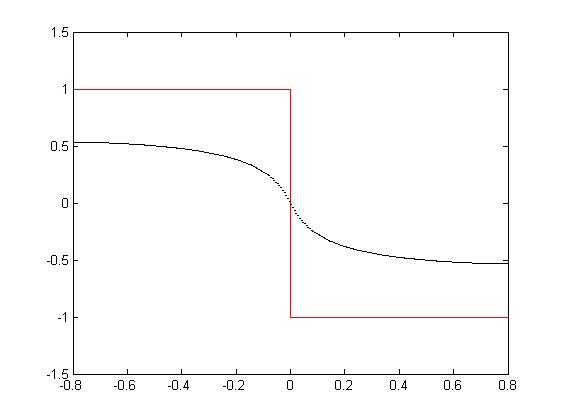}}
\caption{Initial data (piecewise linear) and solutions of the pure fractional equation ($\lambda=0.5$) with $k=0$ and $\Delta x=1/160$.
}\label{pf}
\end{figure}

\example Let us consider the pure fractional equation $\partial_{t}u=g[u]$. From e.g.~\cite{Levy}, it follows that the solution of this equation is given by
the convolution product $u(x,t)=(K*u_{0})(x,t)$, where $K$ is the kernel of $g$. Using the properties of the kernel, it can be shown that this equation has a
regularizing effect on the initial datum (see e.g.~\cite{Alibaud/Droniou/Vovelle}); this regularization appears clearly in our numerical experiments presented
in \textsc{Figure} \ref{pf}.

\begin{figure}[t]
\subfigure[$T=0.1$]{
\includegraphics[scale=0.31]{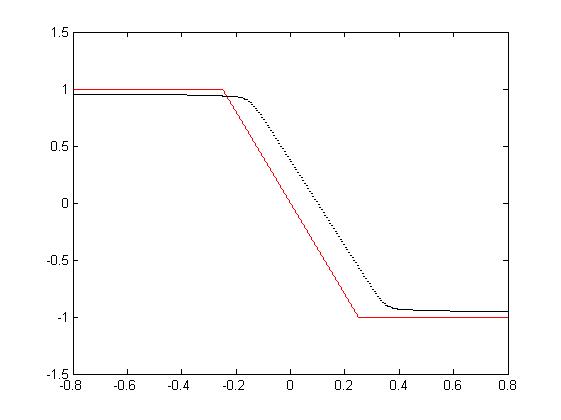}}
\subfigure[$T=0.2$]{
\includegraphics[scale=0.31]{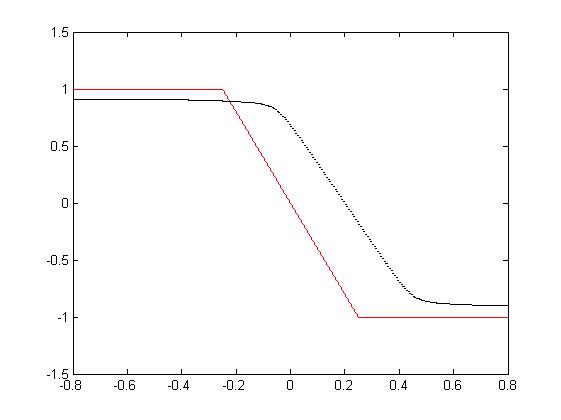}}\\
\subfigure[$T=0.1$]{
\includegraphics[scale=0.31]{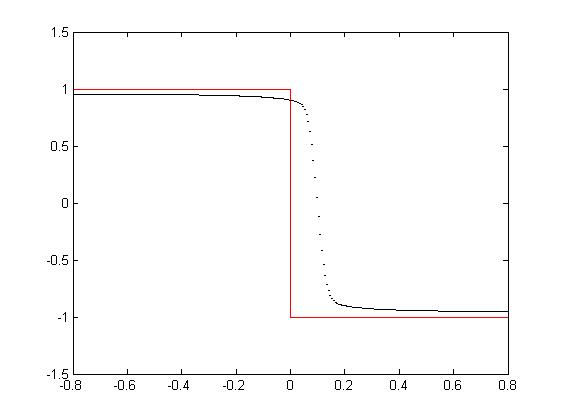}}
\subfigure[$T=0.2$]{
\includegraphics[scale=0.31]{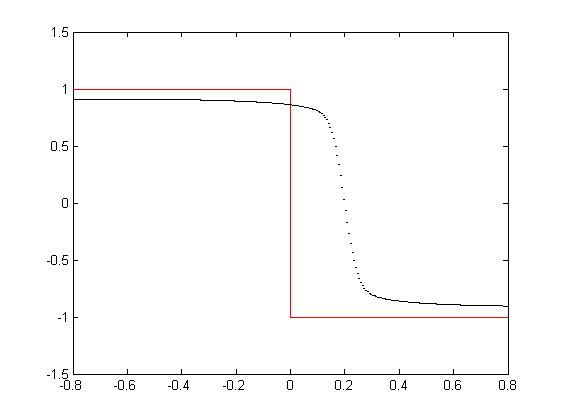}}
\caption{Initial data (piecewise linear) and solutions of the
  fractional transport equation ($\lambda=0.5$) with $k=0$ and $\Delta
  x=1/160$.
}\label{ft}
\end{figure}

\example Let us consider the fractional transport equation $\partial_{t}u+\partial_{x}u=g[u]$. Our numerical results suggest that, as done by
$\partial_{t}u+\partial_{x}u=\partial_{x}^{2}u$, this equation regularizes and transports the initial datum. Our numerical experiments are presented in
\textsc{Figure} \ref{ft}. The numerical flux \eqref{FLUX} has been used.

\begin{figure}[t]
\subfigure[$u_{0}(x)=-\text{sgn}(x)$]{
\includegraphics[scale=0.31]{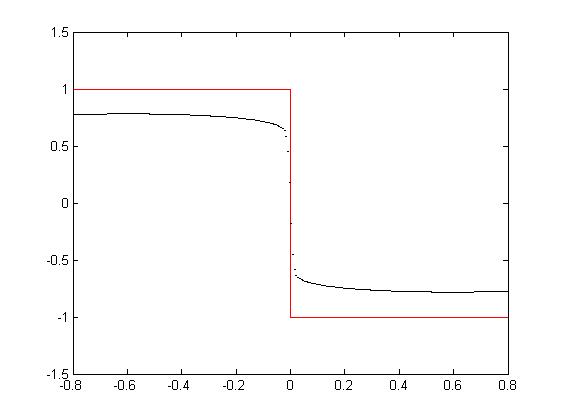}}
\subfigure[$u_{0}(x)=-\text{arctan}(15x)/90$]{
\includegraphics[scale=0.31]{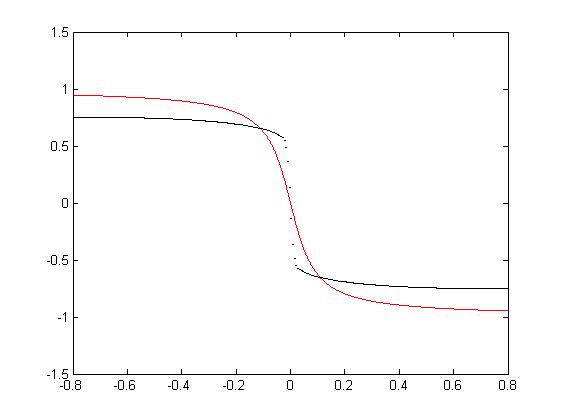}}\\
\subfigure[$u_{0}(x)=\text{sgn}(x)\mathbf{1}_{|x|>1/4}+4x\mathbf{1}_{|x|\leq1/4}$]{
\includegraphics[scale=0.31]{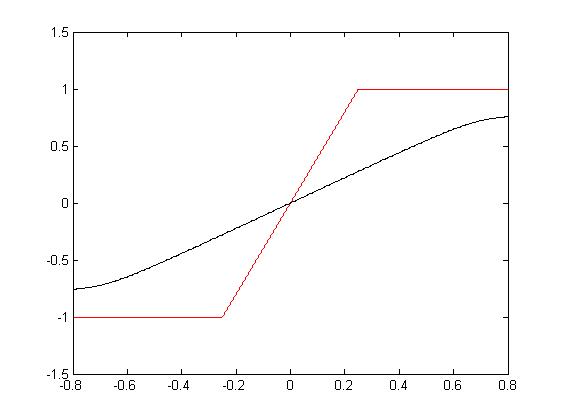}}
\subfigure[$u_{0}(x)=\sin(2\pi x)$]{
\includegraphics[scale=0.31]{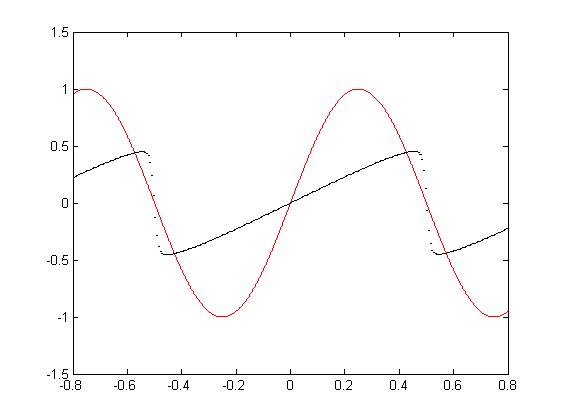}}
\caption{Initial data and solutions of the fractional Burgers'
  equation ($\lambda=0.5$) using $k=0$;
  $T=0.5$ and $\Delta x=1/160$.}\label{fb}
\end{figure}

\begin{figure}[t]
\subfigure[$\lambda=0.1$]{
\includegraphics[scale=0.31]{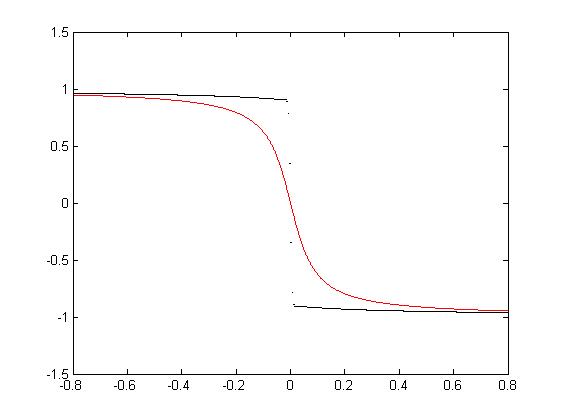}}
\subfigure[$\lambda=0.3$]{
\includegraphics[scale=0.31]{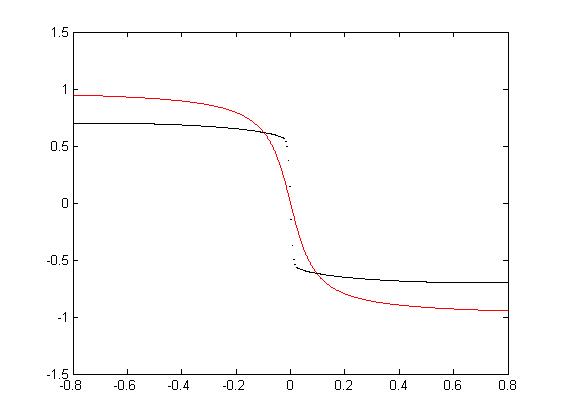}}\\
\subfigure[$\lambda=0.7$]{
\includegraphics[scale=0.31]{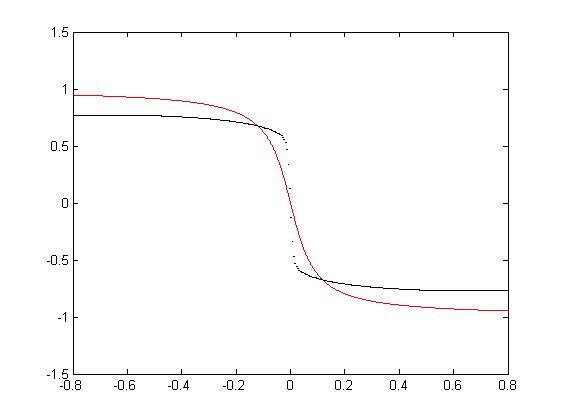}}
\subfigure[$\lambda=0.99$]{
\includegraphics[scale=0.31]{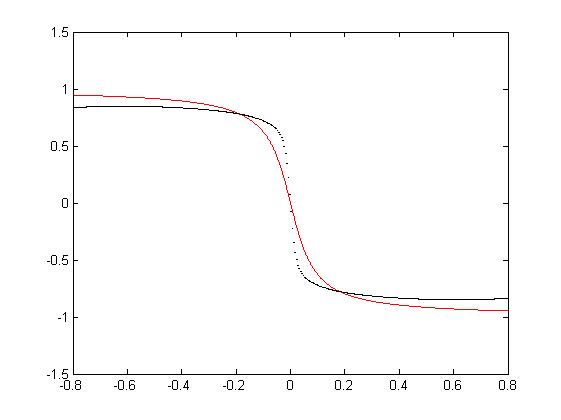}}
\caption{Initial data and solutions of the fractional Burgers'
  equation for different values of $\lambda$ using $k=0$; $T=0.5$,
  $\Delta x=1/200$, and $u_{0}(x)=-\text{arctan}(15x)/90$.
}\label{fbis}
\end{figure}

\begin{figure}[t]
\subfigure[$T=0.1$]{
\includegraphics[scale=0.31]{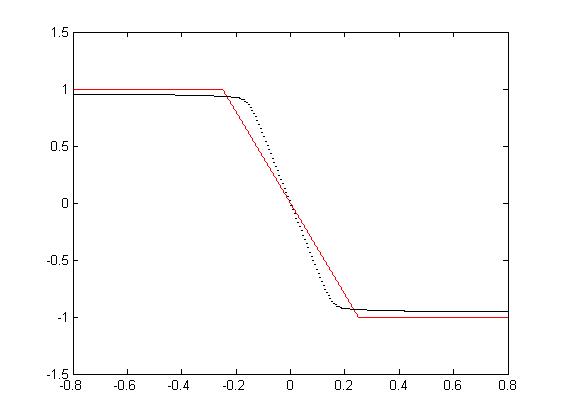}}
\subfigure[$T=0.7$]{
\includegraphics[scale=0.31]{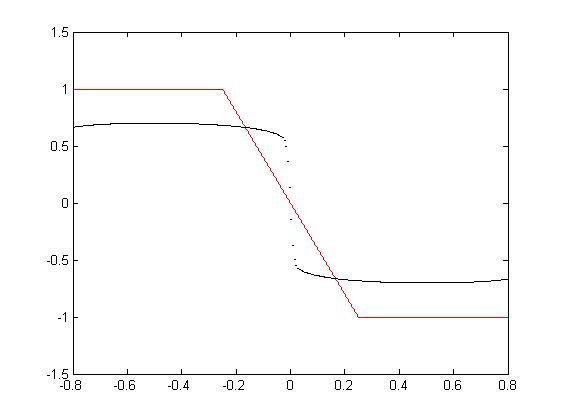}}\\
\subfigure[$T=1.7$]{
\includegraphics[scale=0.31]{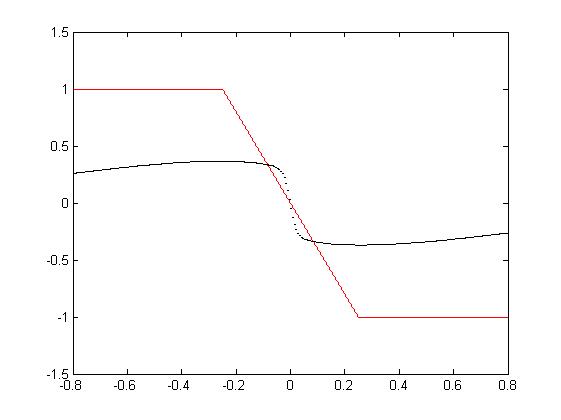}}
\subfigure[$T=2.9$]{
\includegraphics[scale=0.31]{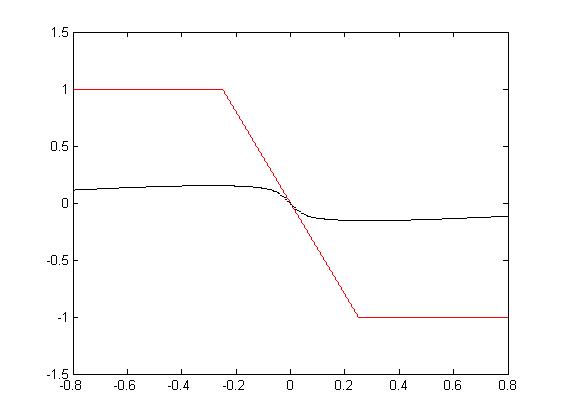}}
\caption{Initial data (piecewise linear) and solutions of the
  fractional Burgers' equation ($\lambda=0.5$) at different times $T$
  using $k=0$; $\Delta x=1/200$.
}\label{fbiss}
\end{figure}

\begin{figure}[t]
\subfigure[$k=0$]{
\includegraphics[scale=0.31]{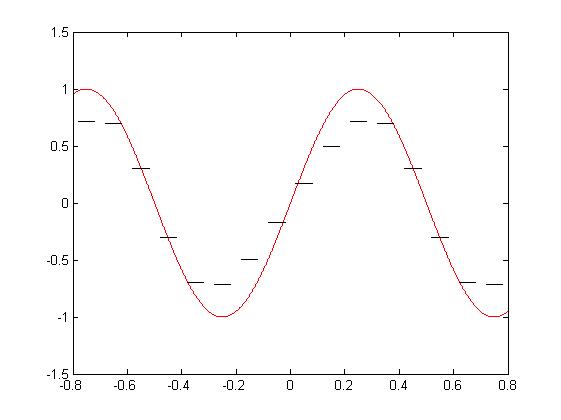}}
\subfigure[$k=1$]{
\includegraphics[scale=0.31]{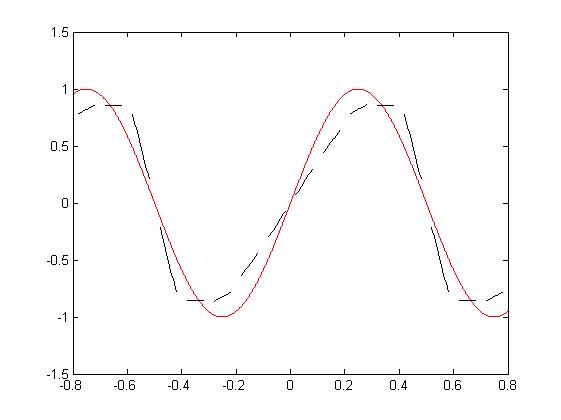}}\\
\subfigure[$k=2$]{
\includegraphics[scale=0.31]{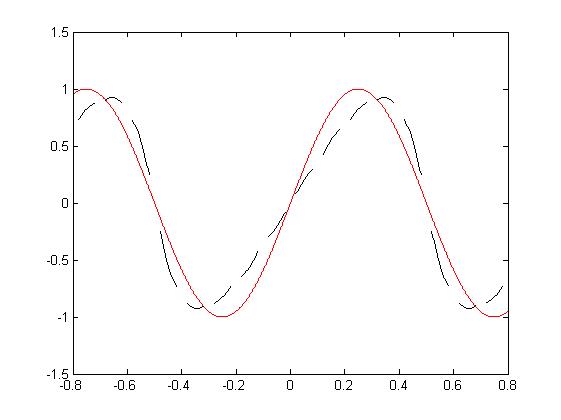}}
\subfigure[Solution computed using $\Delta x=1/640$]{
\includegraphics[scale=0.31]{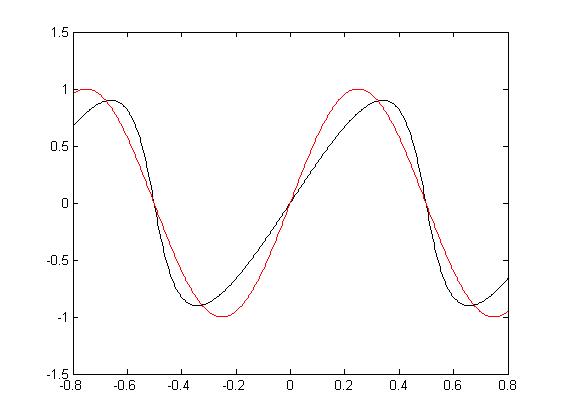}}
\caption{Initial data and solutions of the fractal
  Burgers' equation at $T=1/10$ using different values of $k=0,1,2$;
  $\Delta x=1/10$, and $u_{0}(x)=\sin(2\pi x)$.}\label{seq2}
\end{figure}

\example Let us consider the fractional Burgers' equation $\partial_{t}u+u\partial_{x}u=g[u]$. Our numerical experiments in \textsc{Figure} \ref{fb} confirm
what has been shown by \cite{Alibaud/Droniou/Vovelle,Kiselev/Nazarov/Shterenberg}: this equation does not regularize the initial condition. Discontinuities in
the initial datum can persist in the solution, and shocks can develop from smooth initial data. \textsc{Figure} \ref{fbis} shows how the behavior of the
solution changes with $\lambda$: as $\lambda\rightarrow0$, our numerical solution approaches the solution of the pure Burgers' equation with a source,
$\partial_{t}u+u\partial_{x}u=u$; as $\lambda\rightarrow1$, our numerical solution approaches the smooth solution of the fractional Burgers' equation with
$\lambda=1$ (see \cite{Kiselev/Nazarov/Shterenberg}). \textsc{Figure} \ref{fbiss} clearly shows how a shock can develop and vanish in a finite time.
\textsc{Figure} \ref{seq2} shows how the accuracy improves with $k=0,1,2$. A third order Runge-Kutta (RK3) time discretization and slope limiters
(cf.~\cite{Cockburn}) have been deployed in \textsc{Figure} \ref{seq2}. We have used the Lax-Friedrichs flux
$$F(a,b)=\frac{1}{2}[f(a)+f(b)-c(b-a)],\quad c=\max\{|f'(a)|:|a|\leq\|u_{0}\|_{L^{\infty}(\mathbb{R})}\}.$$ Let us note that the above numerical flux does not
fulfil assumption \emph{A1}. However, this assumption can be replaced with a milder one: it is enough to ask $F(a,b)$ to be Lipschitz continuous on
$\{(a,b):|a|\leq\|u_{0}\|_{L^{\infty}(\mathbb{R})}\text{ and }|b|\leq\|u_{0}\|_{L^{\infty}(\mathbb{R})}\}$.

To give an idea about the speed of convergence of our experiments, we have computed their rate of convergence in \textsc{Table} \ref{table}. We have measured
the error $$E_{\Delta x,p}:=\|\tilde{u}_{\Delta x}(\cdot,T)-\tilde{u}_{e}(\cdot,T)\|_{L^{p}(\mathbb{R})}$$ ($\tilde{u}_{e}$ is the numerical solution which has
been computed using $\Delta x=1/640$), the relative error $$R_{\Delta x,p}:=E_{\Delta x,p}/\|\tilde{u}_{e}(\cdot,T)\|_{L^{p}(\mathbb{R})},$$ and the
approximate rate of convergence $$\alpha_{\Delta x,p}:=(\log E_{\Delta x,p}-\log E_{\Delta x/2,p})/\log 2.$$

We expected to see numerical convergence of order $1/2$ for $k=0$ and numerical convergence of order $3/2$ for $k=1$ (i.e, high-order convergence). The values
$\alpha_{\Delta x,1}$ roughly suggest $1/2$ convergence while the values $\alpha_{\Delta x,2}$ do not reach the expected rate $3/2$. This could be due to our
way or reducing the problem from a nonlocal to a local one (cf. Remark \ref{rem}).

\begin{table}
\caption{$k=0$ (left) as in \textsc{Figure} \ref{fb} $(c)$ and $k=1$ (right) as in \textsc{Figure} \ref{seq2} $(b)$.}\label{table}
\begin{tabular}{c| c c c| c c c}
$\mathbf{\Delta x}$ & $\mathbf{E_{\Delta x,1}}$ & $\mathbf{R_{\Delta x,1}}$ & $\mathbf{\alpha_{\Delta x,1}}$
&  $\mathbf{E_{\Delta x,2}}$ &  $\mathbf{R_{\Delta x,2}}$   &   $\mathbf{\alpha_{\Delta x,2}}$\\
    \hline
$\mathbf{1/10}$    & 0.1990 & 0.1109 & 0.5726 & 0.4580 & 0.3765 & 1.0714\\
$\mathbf{1/20}$    & 0.1338 & 0.0746 & 0.4711 & 0.2180 & 0.1792 & 1.2024\\
$\mathbf{1/40}$    & 0.0965 & 0.0538 & 0.3964 & 0.0947 & 0.0779 & 1.1717\\
$\mathbf{1/80}$    & 0.0734 & 0.0409 & 0.4399 & 0.0421 & 0.0346 & 1.0881\\
$\mathbf{1/160}$   & 0.0541 & 0.0301 & 0.7235 & 0.0198 & 0.0163 & -\\
$\mathbf{1/320}$   & 0.0327 & 0.0183 & -      & -      & -      & -\\
    \hline
\end{tabular}
\end{table}


\appendix

\section{Technical lemmas}
\begin{lemma}\label{00012}
Let $\varphi\in L^{1}(\mathbb{R})\cap BV(\mathbb{R})$. Then, there exists $C>0$ such that
\begin{equation*}
\|g[\varphi]\|_{L^{1}(\mathbb{R})}\leq c_{\lambda}\int_{\mathbb{R}}\int_{|z|>0}\frac{|\varphi(x+z)-\varphi(x)|}{|z|^{1+\lambda}}dzdx\leq c_\lambda
C\|\varphi\|_{L^{1}(\mathbb{R})}^{1-\lambda}|\varphi|_{BV(\mathbb{R})}^{\lambda}.
\end{equation*}
\end{lemma}
\begin{proof}
For all $\epsilon>0$,
\begin{align*}
&\int_{|z|<\epsilon}\int_{\mathbb{R}}\frac{|\varphi(x+z)-\varphi(x)|}{|z|^{1+\lambda}}dxdz
\leq\epsilon^{1-\lambda}|\varphi|_{BV(\mathbb{R})}\int_{|z|<1}\frac{1}{|z|^{\lambda}}dz,\\
&\int_{|z|>\epsilon}\int_{\mathbb{R}}\frac{|\varphi(x+z)-\varphi(x)|}{|z|^{1+\lambda}}dxdz \leq
2\epsilon^{-\lambda}\|\varphi\|_{L^{1}(\mathbb{R})}\int_{|z|>1}\frac{1}{|z|^{1+\lambda}}dz.
\end{align*}
Set $\epsilon=\frac{\|\varphi\|_{L^{1}(\mathbb{R})}}{|\varphi|_{BV(\mathbb{R})}}$ to conclude.
\end{proof}

\begin{lemma}\label{000}
Let $\varphi,\phi\in L^{1}(\mathbb{R})\cap BV(\mathbb{R})$. Then
\begin{gather*}
\int_{\mathbb{R}}\varphi g[\phi]dx=\int_{\mathbb{R}}g[\varphi]\phi dx\\
\intertext{and, in particular,}
\int_{\mathbb{R}}\varphi g[\varphi]dx=-\frac{c_\lambda}{2}\int_{\mathbb{R}}\int_{\mathbb{R}}\frac{(\varphi(z)-\varphi(x))^{2}}{|z-x|^{1+\lambda}}dzdx.
\end{gather*}
\end{lemma}
\begin{proof} By Lemma \ref{00012} and the fact that $BV(\mathbb R)\subset
  L^\infty(\mathbb R)$,
\begin{equation*}
\|\varphi g[\phi]\|_{L^1(\mathbb{R})}
\leq c_\lambda C\|\phi\|^{1-\lambda}_{L^1(\mathbb{R})}|\phi|^{\lambda}_{BV(\mathbb{R})}\|\varphi\|_{L^\infty(\mathbb{R})}<\infty,
\end{equation*}
and, then, Fubini's theorem can be used to obtain
\begin{align*}
\int_{\mathbb{R}}\varphi(x)g[\phi(x)]dx=\frac{1}{2}\int_{\mathbb{R}}\int_{\mathbb{R}}\frac{(\phi(x)-\phi(z))(\varphi(z)-\varphi(x))}{|z-x|^{1+\lambda}}dzdx
=\int_{\mathbb{R}}g[\varphi(x)]\phi(x)dx.
\end{align*}
\end{proof}

\begin{corollary}\label{023}
Lemma \ref{000} holds true for all $\varphi,\phi\in H^{\lambda/2}(\mathbb{R})$.
\end{corollary}
\begin{proof}
  Lemma \ref{000} holds true, in particular, for all
  $\varphi_{n},\phi_{n}$ step functions with compact support,
\begin{equation}\label{dens}
\int_{\mathbb{R}}\varphi_{n}(x)g[\phi_{n}(x)]dx=\int_{\mathbb{R}}g[\varphi_{n}(x)]\phi_{n}(x)dx.
\end{equation}
Let us choose, by density, $\varphi_{n},\phi_{n}\rightarrow\varphi,\phi$ in $H^{\lambda/2}(\mathbb{R})$, and recall the following definition of the
$H^{\lambda/2}$-norm (cf.~\cite[Chapter 6]{Folland}):
\begin{align}\label{foll}
\|\varphi\|_{H^{\lambda/2}(\mathbb{R})}^{2}:=\int_{\mathbb{R}}(1+\xi^{2})^{\lambda/2}\hat{\varphi}^{2}(\xi)d\xi.
\end{align}
Note that, using \eqref{foll} and \eqref{fourier},
\begin{align*}
\|g[\varphi_{n}]-g[\varphi]\|_{H^{-\lambda/2}(\mathbb{R})}&=\int_{\mathbb{R}}(1+\xi^{2})^{-\lambda/2}\xi^{2\lambda}[\hat{\varphi}_{n}(\xi)-\hat{\varphi}(\xi)]^{2}d\xi\\
&\leq\int_{\mathbb{R}}(1+\xi^{2})^{\lambda/2}[\hat{\varphi}_{n}(\xi)-\hat{\varphi}(\xi)]^{2}d\xi=\|\varphi_{n}-\varphi\|_{H^{\lambda/2}(\mathbb{R})}
\end{align*}
since $(1+\xi^{2})^{-\lambda/2}\xi^{2\lambda}\leq(1+\xi^{2})^{\lambda/2}$ for all $\xi\in\mathbb{R}$ (indeed, call $\xi^2=x$, and multiply both sides by
$(1+x)^{1-\lambda/2}$ to get $x^\lambda\leq(1+x)^\lambda$ which holds true for all $x\geq0$). Thus, since $g[\varphi_{n}],g[\phi_{n}]\rightarrow g[\varphi],g[\phi]$
in $H^{-\lambda/2}(\mathbb{R})$ whenever $\varphi_n,\phi_n\rightarrow\varphi,\phi$ in $H^{\lambda/2}(\mathbb{R})$, equality \eqref{dens} holds true also in the
limit $n\rightarrow\infty$.
\end{proof}

\begin{lemma}\label{lem:new}
If $\phi\in V^k\cap L^{2}(\mathbb{R})$, then $\phi\in H^{\lambda/2}(\mathbb{R})$ and, for some constant $c>0$,
$$\|\phi\|^2_{H^{\lambda/2}(\mathbb{R})}\leq\frac{c}{\Delta x}\|\phi\|^2_{L^{2}(\mathbb{R})}.$$
\end{lemma}
\begin{proof}
Let us choose a function $\phi\in V^k\cap L^2(\mathbb{R})$,
$\phi(x)=\sum_{i\in\mathbb{Z}}\sum_{p=0}^k c_{p,i}\varphi_{p,i}(x)$,
and let $\phi_r'$ be the regular part of its derivative,
$$\phi_r'(x)=\sum_{i\in\mathbb{Z}}\sum_{p=0}^k c_{p,i}\frac{d}{dx}\varphi_{p,i}(x).$$ In this case we may define the quadratic variation of $\phi$ as
\begin{equation*}
\begin{split}
|\phi|^2_{QV(\mathbb{R})}=\sum_{i\in\mathbb{Z}}[\phi(x_i^+)-\phi(x_i^-)]^2.
\end{split}
\end{equation*}
First of all, we prove that both $\|\phi'_r\|_{L^2(\mathbb{R})},|\phi|_{QV(\mathbb{R})}$ are finite since $\phi\in L^2(\mathbb{R})$. Indeed, by
orthogonality of the Legendre polynomials,
$$\|\phi\|^2_{L^2(\mathbb{R})}=\sum_{i\in\mathbb{Z}}\sum_{p=0}^k
\frac{\Delta x}{2p+1}c^2_{p,i}\quad\text{and hence}\quad
\sum_{i\in\mathbb{Z}}\sum_{p=0}^k c^2_{p,i}\leq\frac{2k+1}{\Delta x}\|\phi\|^2_{L^2(\mathbb{R})}.$$ As a consequence $\phi'_r\in L^2(\mathbb{R})$ since for each Legendre
polynomial $\varphi_{p,i}$,
$\frac{d}{dx}\varphi_{p,i}=\sigma_p\varphi_{p-1,i}$ for some constant $\sigma_p$. Moreover,
\begin{equation*}
\begin{split}
|\phi|^2_{QV(\mathbb{R})}\leq \frac{(k+1)(2k+1)}{\Delta x}\|\phi\|^2_{L^2(\mathbb{R})}
\end{split}
\end{equation*}
since $|\phi|^2_{QV(\mathbb{R})}\leq 2\sum_{i\in\mathbb{Z}}\phi^2(x_i^+)+2\sum_{i\in\mathbb{Z}}\phi^2(x_i^-)$ and (remember that $\varphi_{p,i}(x_{i+1}^{-})=1$
while $\varphi_{p,i}(x_{i}^{+})=(-1)^{p}$)
\begin{equation*}
\begin{split}
\sum_{i\in\mathbb{Z}}\phi^2(x_i^\pm)&=\sum_{i\in\mathbb{Z}}\left(\sum_{p=0}^k c_{p,i}\varphi_{p,i}(x_i^\pm)\right)^2\\
&\leq (k+1)\sum_{i\in\mathbb{Z}}\sum_{p=0}^k c_{p,i}^2\leq\frac{(k+1)(2k+1)}{\Delta x}\|\phi\|^2_{L^2(\mathbb{R})}.
\end{split}
\end{equation*}
Next, we prove that there exists a constant $c>0$ such that, for a.e.~$|z|<1$,
\begin{equation}\label{kaka}
\begin{split}
\int_{\mathbb{R}}[\phi(x+z)-\phi(x)]^2dx\leq c\left(|z||\phi|^2_{QV(\mathbb{R})}+|z|^2\|\phi'_r\|^2_{L^2(\mathbb{R})}\right).
\end{split}
\end{equation}
Note that
\begin{equation*}
\begin{split}
\int_{\mathbb{R}}[\phi(x+z)-\phi(x)]^2dx&=\sum_{i\in\mathbb{Z}}\int_{i|z|}^{(i+1)|z|}[\phi(x+z)-\phi(x)]^2dx\\
&=\int_{0}^{|z|}\sum_{i\in\mathbb{Z}}[\phi(x+(i+1)|z|)-\phi(x+i|z|)]^2dx.
\end{split}
\end{equation*}
By appropriately adding and subtracting the values $\phi(x_i^\pm)$, $i\in\mathbb{Z}$, the right-hand side of the above expression is less than or equal to
\begin{equation}\label{vvv}
\begin{split}
3\int_{0}^{|z|}\sum_{i\in\mathbb{Z}}[\phi(x_i^+)-\phi(x_i^-)]^2dx+3\int_{0}^{|z|}\sum_{i\in\mathbb{Z}}\sum_{j=0}^{J_i-1}[\phi(z_{j}^i)-\phi(z_{j+1}^i)]^2dx,
\end{split}
\end{equation}
where the $J_i+1$ points $x_i^+=z_0^i\leq\ldots\leq z_{J_i}^i=x_{i+1}^-$ lie inside the interval $I_i$ (these points can vary from
interval to interval depending on the value of $\Delta x$. E.g. if
$|z|\ll \Delta x$, each interval $I_i$ contains more than two points,
while if $|z|\gg \Delta x$, some intervals contain just the end-points
$x_i^+=z_0^i$ and $z_{J_i}^i=x_{i+1}^-$. We can control the first
term in \eqref{vvv} 
thanks to the bound on the quadratic variation of $\phi$ while, since inside each interval $I_i$ the function $\phi$ is smooth, we can use the Taylor's formula to
rewrite the second term as
$$\sum_{i\in\mathbb{Z}}\sum_{j=0}^{J_i-1}[\phi(z_{j}^i)-\phi(z_{j+1}^i)]^2\leq|z|\sum_{i\in\mathbb{Z}}\sum_{j=0}^{J_i-1}[\phi'_r(y_{j}^i)]^2(z_{j+1}^i-z_j^i),
\text{ where $z_j^i\leq y_{j}^i\leq z_{j+1}^i$.}$$
The right-hand side of the above inequality contains
a Riemann sum approximation of the $L^2$-norm of the function $\phi'_r\in V^k\cap L^2(\mathbb{R})$
and is therefore finite. Hence inequality \eqref{kaka} has been established, and we are now
ready to conclude the proof. The seminorm
\begin{equation*}
\begin{split}
&|\phi|^2_{H^{\lambda/2}(\mathbb{R})}\\
&=\int_{|z|>1}\int_{\mathbb{R}}\frac{[\phi(x+z)-\phi(x)]^2}{|z|^{1+\lambda}}dxdz+\int_{|z|<1}\int_{\mathbb{R}}\frac{[\phi(x+z)-\phi(x)]^2}{|z|^{1+\lambda}}dxdz:=J_1+J_2
\end{split}
\end{equation*}
is finite since
\begin{align*}
J_1\leq 4\|\phi\|^2_{L^{2}(\mathbb{R})}\int_{|z|>1}\frac{dz}{|z|^{1+\lambda}}<\infty,
\end{align*}
and, thanks to \eqref{kaka},
\begin{align*}
J_2\leq c\left(|\phi|^2_{QV(\mathbb{R})}\int_{|z|<1}\frac{dz}{|z|^{\lambda}}+\|\phi'_r\|^2_{L^2(\mathbb{R})}\int_{|z|<1}\frac{dz}{|z|^{\lambda-1}}\right)<\infty.
\end{align*}
\end{proof}

\begin{lemma}\label{00012bis}
Let $u\in H^{k+1}(\mathbb{R})$ and $\mathbf{u}$ be its $L^2$-projection into $V^k$, then
there exists a constant $c_k>0$ such that $$\|u-\mathbf{u}\|_{H^{\lambda/2}(\mathbb{R})}^{2}\leq c_{k}\|u\|^2_{H^{k+1}(\mathbb{R})}\Delta
x^{2k+2-\lambda}.$$
\end{lemma}
\begin{proof}
Let us call $v=u-\mathbf{u}$, and remember that (cf.~\cite[Section
4.4]{Brenner/Scott}), for some constant $c_k>0$ and all intervals
$I_i=(i\Delta x,(i+1)\Delta x)$,
\begin{align*}
&\|v\|_{L^{2}(I_i)}\leq c_{k}\|u\|_{H^{k+1}(I_i)}\Delta x^{k+1},\\
&\|v\|_{L^{\infty}(I_i)}\leq c_{k}\|u\|_{H^{k+1}(I_i)}\Delta x^{k+\frac{1}{2}},\\
&\|v\|_{H^{1}(I_i)}\leq c_{k}\|u\|_{H^{k+1}(I_i)}\Delta x^{k}.
\end{align*}
First of all, let us bound from above the $H^{\lambda/2}$-norm of $v$ as
\begin{equation*}
\begin{split}
\|v\|^2_{H^{\lambda/2}(\mathbb{R})}&\leq\|v\|^2_{L^{2}(\mathbb{R})}+\sum_{i\in\mathbb{Z}}\int_{I_i}\int_{I_i}\frac{[v(z)-v(x)]^2}{|z-x|^{1+\lambda}}dzdx\\
&\qquad\qquad\qquad+2\sum_{i\in\mathbb{Z}}\int_{I_i}\int_{I_{i+1}}\frac{[v(z)-v(x)]^2}{|z-x|^{1+\lambda}}dzdx\\
&\qquad\qquad\qquad+\int_{\mathbb{R}}\int_{|z-x|>\Delta x}\frac{[v(z)-v(x)]^2}{|z-x|^{1+\lambda}}dzdx\\
&:=J_1+J_2+J_3+J_4.
\end{split}
\end{equation*}
Note that
\begin{equation*}
\begin{split}
J_1&=\sum_{i\in\mathbb{Z}}\|v\|^2_{L^{2}(I_i)}\\
&\leq c_k\Delta x^{2k+2}\sum_{i\in\mathbb{Z}}\|u\|^2_{H^{k+1}(I_i)}=c_k\|u\|^2_{H^{k+1}(\mathbb{R})}\Delta x^{2k+2}.
\end{split}
\end{equation*}
We now prove the remaining $J_i$ ($i=2,3,4$) to be of order $\Delta x^{2k+2-\lambda}$. First note that since $v$ is smooth on each interval $I_i$, the
fundamental theorem of calculus followed by Jensen's inequality yield
\begin{equation*}
\begin{split}
\int_{I_i}\int_{I_i}\frac{[v(z)-v(x)]^2}{|z-x|^{1+\lambda}}dzdx&=\int_{I_i}\int_{I_i}\frac{1}{|z-x|^{1+\lambda}}\left(\int_x^z\frac{d}{ds}v(s)ds\right)^2dzdx\\
&\leq\int_{I_i}\int_{I_i}\frac{|z-x|}{|z-x|^{1+\lambda}}\int_{I_i}\left(\frac{d}{ds}v(s)\right)^2dsdzdx\\
&\leq\|v\|^2_{H^1(I_i)}\int_{x_i}^{x_{i+1}}\int_{x_i}^{x_{i+1}}\frac{dzdx}{|z-x|^{\lambda}}\\
&\leq\|v\|^2_{H^1(I_i)}\int_{x_i}^{x_{i+1}}\int_{-\Delta x}^{\Delta x}\frac{dsdx}{|s|^{\lambda}}\\
&\leq\Delta x^{2-\lambda}\|v\|^2_{H^1(I_i)}\int_{-1}^1\frac{ds}{|s|^{\lambda}}.
\end{split}
\end{equation*}
Thus $I_2=c\Delta x^{2-\lambda}\sum_{i\in\mathbb{Z}}\|v\|^2_{H^1(I_i)}\leq c_k\|u\|^2_{H^{k+1}(\mathbb{R})}\Delta x^{2k+2-\lambda}$. Next, we note that
\begin{equation*}
\begin{split}
&\int_{I_i}\int_{I_{i+1}}\frac{[v(z)-v(x)]^2}{|z-x|^{1+\lambda}}dzdx\\
&\leq2\int_{I_i}\int_{I_{i+1}}\frac{[v(z)]^2}{|z-x|^{1+\lambda}}dzdx +2\int_{I_i}\int_{I_{i+1}}\frac{[v(x)]^2}{|z-x|^{1+\lambda}}dzdx.
\end{split}
\end{equation*}
Let us show how to estimate the first term on the right-hand side of the above inequality. Analogous ideas can be used for the second one.
\begin{equation*}
\begin{split}
\int_{I_i}\int_{I_{i+1}}\frac{[v(z)]^2}{|z-x|^{1+\lambda}}dzdx&\leq\|v\|^2_{L^\infty(I_{i+1})}\int_{x_i}^{x_{i+1}}\int_{x_{i+1}}^{x_{i+2}}\frac{dzdx}{(z-x)^{1+\lambda}}\\
&\leq\|v\|^2_{L^\infty(I_{i+1})}\int_{x_i}^{x_{i+1}}\int_{x_{i+1}-x}^\infty\frac{dsdx}{s^{1+\lambda}}\\
&=\|v\|^2_{L^\infty(I_{i+1})}\int_{1}^\infty\frac{dr}{r^{1+\lambda}}\int_{x_i}^{x_{i+1}}\frac{dx}{(x_{i+1}-x)^{\lambda}}\\
&\leq\|v\|^2_{L^\infty(I_{i+1})}\int_{1}^\infty\frac{dr}{r^{1+\lambda}}\int_0^{\Delta x}\frac{dy}{y^{\lambda}}\\
&=\|v\|^2_{L^\infty(I_{i+1})}\int_{1}^\infty\frac{dr}{r^{1+\lambda}}\int_{0}^1\frac{dy}{y^{\lambda}}\
\Delta x^{1-\lambda}.\\
\end{split}
\end{equation*}
Thus $I_3=c\Delta x^{1-\lambda}\sum_{i\in\mathbb{Z}}\|v\|^2_{L^\infty(I_i)}\leq c_k\|u\|^2_{H^{k+1}(\mathbb{R})}\Delta x^{2k+2-\lambda}$. Finally,
\begin{equation*}
\begin{split}
\int_{\mathbb{R}}\int_{|z-x|>\Delta x}\frac{[v(z)-v(x)]^2}{|z-x|^{1+\lambda}}dzdx&=\int_{\mathbb{R}}\int_{|s|>\Delta
x}\frac{[v(x+s)-v(x)]^2}{|s|^{1+\lambda}}dsdx\\
&\leq 4\|v\|^2_{L^2(\mathbb{R})}\Delta x^{-\lambda}\int_{|s|>1}\frac{ds}{|s|^{1+\lambda}},
\end{split}
\end{equation*}
and $I_4=c\Delta x^{-\lambda}\sum_{i\in\mathbb{Z}}\|v\|^2_{L^2(I_i)}\leq c_k\|u\|^2_{H^{k+1}(\mathbb{R})}\Delta x^{2k+2-\lambda}$.
\end{proof}

\begin{lemma}\label{cor:s}
Let $u\in V^k\cap L^2(\mathbb{R})$, $a_{p,i}=\int_{I_i}g[u]\varphi_{p,i}$ and $$\gamma_u(x)=\sum_{i\in\mathbb{Z}}\sum_{p=0}^ka_{p,i}\varphi_{p,i}(x).$$ Then,
$\|\gamma_u\|_{L^2(\mathbb{R})}\leq c\|u\|_{L^2(\mathbb{R})}$ for some constant $c>0$.
\end{lemma}
\begin{proof}
Let us introduce the compactly supported function $v_M\in V^k\cap L^2(\mathbb{R})$ $$v_M(x)=\sum_{|i|\leq M}\sum_{p=0}^k a_{p,i}\varphi_{p,i}(x).$$ Note that,
since $a_{p,i}=\int_{I_i}g[u]\varphi_{p,i}$,
\begin{equation*}
\begin{split}
\sum_{|i|\leq M}\sum_{p=0}^ka_{p,i}^2&=\sum_{|i|\leq M}\sum_{p=0}^ka_{p,i}\int_{I_i}g[u]\varphi_{p,i}=\int_{\mathbb{R}}g[u]v_M.
\end{split}
\end{equation*}
By Lemma \ref{lem:new}, the pairing $\int_{\mathbb{R}}g[u]v_M$ is less than or equal to
\begin{equation*}
\begin{split}
\|u\|_{H^{\lambda/2}(\mathbb{R})}\|v_M\|_{H^{\lambda/2}(\mathbb{R})}
\leq c\|u\|_{L^2(\mathbb{R})}\|v_M\|_{L^2(\mathbb{R})}\leq c\|u\|_{L^2(\mathbb{R})}\left(\sum_{|i|\leq M}\sum_{p=0}^k a_{p,i}^2\right)^{\frac{1}{2}}
\end{split}
\end{equation*}
for some $c>0$. Hence, $\sum_{|i|\leq M}\sum_{p=0}^ka_{p,i}^2\leq c\|u\|^2_{L^2(\mathbb{R})}$ and, in the limit $M\rightarrow\infty$,
$$\|\gamma_u\|^2_{L^2(\mathbb{R})}=\sum_{i\in\mathbb{Z}}\sum_{p=0}^ka_{p,i}^2\leq c\|u\|^2_{L^2(\mathbb{R})}.$$
\end{proof}

\section{Proof of Proposition \ref{007}}
Since $G_{j}^{i}:=\int_{\mathbb{R}}\mathbf{1}_{I_{i}}(x)g[\mathbf{1}_{I_{j}}(x)]dx$, Lemma \ref{000} returns
\begin{align*}
G_{j}^{i}=\int_{\mathbb{R}}\mathbf{1}_{I_{i}}(x)g[\mathbf{1}_{I_{j}}(x)]dx=\int_{\mathbb{R}}\mathbf{1}_{I_{j}}(x)g[\mathbf{1}_{I_{i}}(x)]dx=G_{i}^{j}.
\end{align*}
Thus, by Lemma \ref{00012}, $\sum_{j\in\mathbb{Z}}|G_{j}^{i}|\leq\int_{\mathbb{R}}|g[\mathbf{1}_{I_{i}}(x)]|dx<\infty$ and, by symmetry,
\begin{align*}
\sum_{j\in\mathbb{Z}}G_{j}^{i}=c_{\lambda}\int_{\mathbb{R}}\int_{\mathbb{R}}\frac{\mathbf{1}_{I_{i}}(z)-\mathbf{1}_{I_{i}}(x)}{|z-x|^{1+\lambda}}dzdx=0.
\end{align*}
All diagonal elements are equal and negative. Indeed,
\begin{equation*}
G_{i}^{i}=c_{\lambda}\int_{I_{i}}\int_{|z|>0}\frac{\mathbf{1}_{I_{i}}(x+z)-\mathbf{1}_{I_{i}}(x)}{|z|^{1+\lambda}}dzdx
=c_{\lambda}\int_{|z|>0}\frac{\xi(z)}{|z|^{1+\lambda}}dz,
\end{equation*}
where
\begin{equation*}
\xi(z)=\left\{
\begin{array}{lll}
-|z|&z\in(-\Delta x,\Delta x)\\
-\Delta x&otherwise.
\end{array}
\right.
\end{equation*}
Thus, $G_{i}^{i}=-c_{\lambda}(\int_{|z|<1}\frac{1}{|z|^{\lambda}}dz+\int_{|z|>1}\frac{1}{|z|^{1+\lambda}}dz)\Delta x^{1-\lambda}$. All elements outside the
diagonal are positive. Moreover, $G_{j+1}^{i+1}=G_{j}^{i}$ for all $(i,j)\in\mathbb{Z}\times\mathbb{Z}$ since, if $i\neq j$,
\begin{equation*}
G_{j}^{i}=c_{\lambda}\int_{I_{i}}\int_{|z|>0}\frac{\mathbf{1}_{I_{j}}(x+z)}{|z|^{1+\lambda}}dzdx.
\end{equation*}

\section{Proof of Theorem \ref{025} for the implicit-explicit method}
Let us consider the problem
\begin{equation}\label{aux}
\begin{split}
v_i-\Delta tg\langle v\rangle_{i}=h_i,\text{ $i\in\mathbb{Z}$ and $h\in l^\infty(\mathbb{Z})\cap l^1(\mathbb{Z})$.}
\end{split}
\end{equation}
One can proceed as done by Droniou for nonlocal operators satisfying all the assumptions listed in \cite{Droniou} (cf.~also \cite{Cifani} for a detailed proof
for the operator $g\langle \cdot\rangle$) to prove the existence of a solution $v\in l^\infty(\mathbb{Z})\cap l^1(\mathbb{Z})$ of problem \eqref{aux}.
Moreover,
\begin{align}
\inf_{i\in\mathbb{Z}}h_i\leq\inf_{i\in\mathbb{Z}}v_i&\leq\sup_{i\in\mathbb{Z}}v_i\leq\sup_{i\in\mathbb{Z}}h_i,\label{111}\\
\sum_{i\in\mathbb{Z}}|v_i|&\leq\sum_{i\in\mathbb{Z}}|h_i|\label{222}.
\end{align}
Note that \eqref{111} ensures uniqueness for problem \eqref{aux}. Our plan is to rewrite the implicit-explicit method \eqref{implicit} in the form
\eqref{aux}, and use \eqref{111}-\eqref{222} to prove Theorem \ref{045}. We start by rewriting \eqref{implicit} in
"linearized" form,
\begin{equation*}
\begin{split}
U^{n+1}_i-\frac{\Delta t}{\Delta x}\sum_{j\in\mathbb{Z}}G_j^iU^{n+1}_j&=U_{i}^{n}-\Delta tD_-F(U_{i}^{n},U_{i+1}^{n})\\
&=a_iU_{i+1}^{n}+(1-a_i-b_i)U_{i}^{n}+b_iU_{i-1}^{n},
\end{split}
\end{equation*}
where
\begin{equation*}
\begin{split}
a_i^n=-\frac{\Delta t}{\Delta x}\frac{F(U^n_i,U^n_{i+1})-F(U^n_i,U^n_i)}{U^n_{i+1}-U^n_{i}}\text{ and }
b_i^n=\frac{\Delta t}{\Delta
x}\frac{F(U^n_{i-1},U^n_{i})-F(U^n_i,U^n_i)}{U^n_{i-1}-U^n_{i}}.
\end{split}
\end{equation*}
(the above coefficients are equal to zero when the denominators are equal to zero). By the CFL condition \eqref{CFLim} and the Lipschitz regularity of $F$, it
follows that $a_i,b_i,1-a_i-b_i$ are bounded and positive. Thus, the implicit-explicit method \eqref{implicit} reduces to \eqref{aux} if we choose
$v_i:=U^{n+1}_i$ and $h_i:=a_iU_{i+1}^{n}+(1-a_i-b_i)U_{i}^{n}+b_iU_{i-1}^{n}$.

We are now ready to prove Theorem \ref{025}. Items \emph{i} and \emph{ii} are easy consequences of \eqref{111} and \eqref{222}. To prove item \emph{iii},
we call $V_i^n=U^n_{i+1}-U^n_i$, and, by using the implicit-explicit method \eqref{implicit} in linearized form, we obtain
\begin{equation*}
\begin{split}
V^{n+1}_i&+\frac{\Delta t}{\Delta x}\sum_{j\in\mathbb{Z}}G_j^iU^{n+1}_j-\frac{\Delta t}{\Delta x}\sum_{j\in\mathbb{Z}}G_j^{i+1}U^{n+1}_j\\
&=a_{i+1}V_{i+1}^{n}+(1-a_i-b_{i+1})V_{i}^{n}+b_iV_{i-1}^{n}.
\end{split}
\end{equation*}
Note that $\sum_{j\in\mathbb{Z}}G_j^{i+1}U^{n+1}_j=\sum_{j\in\mathbb{Z}}G_{j-1}^{i}U^{n+1}_j=\sum_{j\in\mathbb{Z}}G_{j}^{i}U^{n+1}_{j+1}$ since
$G^{j+1}_{i+1}=G^{j}_{i}$ for all $(i,j)\in\mathbb{Z}\times\mathbb{Z}$, and
\begin{equation*}
\begin{split}
V^{n+1}_i-\frac{\Delta t}{\Delta x}\sum_{j\in\mathbb{Z}}G_j^iV^{n+1}_j=a_{i+1}V_{i+1}^{n}+(1-a_i-b_{i+1})V_{i}^{n}+b_iV_{i-1}^{n},
\end{split}
\end{equation*}
which is of the form \eqref{aux} and, thus, $l^1$-contractive. This proves item \emph{iii}. The proof of item \emph{iv} goes as the one for the fully explicit
method \eqref{explicit}.

\section{Proof of Lemma \ref{045}}
Note that $\Lambda_{\epsilon,\delta}[u,\tilde{u}]\geq0$ by \eqref{BVentropy}, and hence
$\Lambda_{\epsilon,\delta}[\tilde{u},u]\leq\Lambda_{\epsilon,\delta}[\tilde{u},u]+\Lambda_{\epsilon,\delta}[u,\tilde{u}]:=I_1+I_2+I_3+I_4$, where
\begin{equation*}
\begin{split}
I_1:=&\int_{Q_T}\int_{Q_T}\eta_{u(y,s)}(\tilde{u}(x,t))\varphi_t(x,y,t,s)dxdtdyds\\
&+\int_{Q_T}\int_{Q_T}\eta_{\tilde{u}(y,s)}(u(x,t))\varphi_t(x,y,t,s)dxdtdyds,\\
I_2:=&\int_{Q_T}\int_{Q_T}q_{u(y,s)}(\tilde{u}(x,t))\varphi_x(x,y,t,s)dxdtdyds\\
&+\int_{Q_T}\int_{Q_T}q_{\tilde{u}(y,s)}(u(x,t))\varphi_x(x,y,t,s)dxdtdyds,\\
I_3:=&\int_{Q_T}\int_{Q_T}\eta'_{u(y,s)}(\tilde{u}(x,t))g[\tilde{u}(x,t)]\varphi(x,y,t,s)dxdtdyds\\
&+\int_{Q_T}\int_{Q_T}\eta'_{\tilde{u}(y,s)}(u(x,t))g[u(x,t)]\varphi(x,y,t,s)dxdtdyds
\end{split}
\end{equation*}
and
\begin{equation*}
\begin{split}
I_4:=&\int_{Q_T}\int_{\mathbb{R}}\eta_{u(y,s)}(\tilde{u}(x,0))\varphi(x,y,0,s)dxdyds\\
&-\int_{Q_T}\int_{\mathbb{R}}\eta_{u(y,s)}(\tilde{u}(x,T))\varphi(x,y,T,s)dxdyds\\
&+\int_{Q_T}\int_{\mathbb{R}}\eta_{\tilde{u}(y,s)}(u(x,0))\varphi(x,y,0,s)dxdyds\\
&-\int_{Q_T}\int_{\mathbb{R}}\eta_{\tilde{u}(y,s)}(u(x,T))\varphi(x,y,T,s)dxdyds.
\end{split}
\end{equation*}
As shown in \cite[Theorem 3.11]{Holden/Risebro}, $I_1=I_2=0$ while
\begin{equation}\label{holden}
\begin{split}
I_4\leq c(\epsilon+\delta+\Delta x)-\|u(\cdot,T)-\tilde{u}(\cdot,T)\|_{L^1(\mathbb{R})}.
\end{split}
\end{equation}
We now prove that $I_3\leq 0$. Note that, since $g[u]\in{L^{1}(Q_T)}$,
\begin{equation*}
\begin{split}
\int_{Q_T}\int_{Q_T}|\eta'_{\tilde{u}(y,s)}(u(x,t))||g[u(x,t)]|\varphi(x,y,t,s)dxdtdyds<\infty,
\end{split}
\end{equation*}
and we can change the order of integration to obtain
\begin{equation*}
\begin{split}
I_3=&\int_{Q_T}\int_{Q_T}\eta'_{u(y,s)}(\tilde{u}(x,t))g[\tilde{u}(x,t)]\varphi(x,y,t,s)dxdtdyds\\
&+\int_{Q_T}\int_{Q_T}\eta'_{\tilde{u}(x,t)}(u(y,s))g[u(y,s)]\varphi(x,y,t,s)dxdtdyds.
\end{split}
\end{equation*}
Since $\eta'_u(\tilde{u})=-\eta'_{\tilde{u}}(u)$,
\begin{align*}
I_3&=\int_{Q_T}\int_{Q_T}\int_{|z|>0}\text{sgn}(\tilde{u}(x,t)-u(y,s))\varphi(x,y,t,s)\\
&\qquad\qquad\frac{(\tilde{u}(x+z,t)-u(y+z,s))-(\tilde{u}(x,t)-u(y,s))}{|z|^{1+\lambda}}dzdxdtdyds\\
&\leq\int_{Q_T}\int_{Q_T}\int_{|z|>0}\varphi(x,y,t,s)\\
&\qquad\qquad\frac{|\tilde{u}(x+z,t)-u(y+z,s)|-|\tilde{u}(x,t)-u(y,s)|}{|z|^{1+\lambda}}dzdxdtdyds.
\end{align*}
Let us rewrite the right-hand side of the above inequality as a sum of two integrals, and use the change of variables $(z,x,y)\rightarrow(-z,x+z,y+z)$ to obtain
\begin{equation*}
\begin{split}
&\frac{1}{2}\int_{Q_T}\int_{Q_T}\int_{|z|>0}\varphi(x+z,y+z,t,s)\\
&\qquad\qquad\frac{|\tilde u(x,t)-u(y,s)|-|\tilde u(x+z,t)-u(y+z,s)|}{|z|^{1+\lambda}}dzdxdtdyds\\
&+\frac{1}{2}\int_{Q_T}\int_{Q_T}\int_{|z|>0}\varphi(x,y,t,s)\\
&\qquad\qquad\frac{|\tilde u(x+z,t)-u(y+z,s)|-|\tilde u(x,t)-u(y,s)|}{|z|^{1+\lambda}}dzdxdtdyds.
\end{split}
\end{equation*}
By adding up these terms we find that
\begin{equation*}
\begin{split}
&I_3\leq\frac{1}{2}\int_{Q_T}\int_{Q_T}\int_{|z|>0}(\varphi(x+z,y+z,t,s)-\varphi(x,y,t,s))\\
&\qquad\qquad\frac{|\tilde u(x,t)-u(y,s)|-|\tilde u(x+z,t)-u(y+z,s)|}{|z|^{1+\lambda}}dzdxdtdyds,
\end{split}
\end{equation*}
and hence $I_3\leq0$ since $\varphi(x+z,y+z,t,s)=\varphi(x,y,t,s)$. To conclude, let us point out that the following result is needed in
\cite[Theorem 3.11]{Holden/Risebro} to prove \eqref{holden}.

\begin{proposition}
Let $u$ be a BV entropy solution of \Ref{1}. Then, there exists a constant $c>0$
such that $\|u(\cdot,t+\delta)-u(\cdot,t)\|_{L^{1}(\mathbb{R})}\leq
c\delta$.
\end{proposition}
\begin{proof}
Let $0<a<b<T$ and $\mathbf{1}_{[a,b]}^{\epsilon}:\mathbb{R}\rightarrow\mathbb{R}$ be a smooth approximation of $\mathbf{1}_{[a,b]}$. Let us call
$\varphi^{\epsilon}(x,t)=\phi(x)\mathbf{1}_{[a,b]}^{\epsilon}(t)$, where $\phi\in C_{c}^{\infty}(\mathbb{R})$. Thus,
\begin{align*}
\int_{0}^{T}\int_{\mathbb{R}}u\varphi_{t}^{\epsilon}+f(u)\varphi_{x}^{\epsilon}+ug[\varphi^{\epsilon}]dxdt=0
\end{align*}
since $u$ is a BV entropy solution of \Ref{1} and, so, a weak solution (cf.~\cite{Alibaud} for the definition of weak solution). The limit for
$\epsilon\rightarrow0$ is, cf.~\cite[Theorem 7.10]{Holden/Risebro},
\begin{align*}
\int_{\mathbb{R}}\phi(x)[u(x,a)-u(x,b)]dx+\int_{a}^{b}\int_{\mathbb{R}}f(u)\phi_{x}+ug[\phi]dxdt=0
\end{align*}
and
\begin{align*}
\|u(\cdot,b)-u(\cdot,a)&\|_{L^{1}(\mathbb{R})}=\sup_{|\phi|\leq1}\int_{\mathbb{R}}\phi(x)[u(x,b)-u(x,a)]dx\\
&=\sup_{|\phi|\leq1}\Bigg\{-\int_{a}^{b}\int_{\mathbb{R}}f(u)\phi_{x}+ug[\phi]dxdt\Bigg\}\\
&\leq c|u_{0}|_{BV(\mathbb{R})}(b-a)+\sup_{|\phi|\leq1}\Bigg\{-\int_{a}^{b}\int_{\mathbb{R}}ug[\phi]dxdt\Bigg\}.
\end{align*}
To conclude the proof, the following estimate is needed:
\begin{align*}
\sup_{|\phi|\leq1}\Bigg\{-\int_{a}^{b}\int_{\mathbb{R}}ug[\phi]dxdt\Bigg\}=\sup_{|\phi|\leq1}\Bigg\{-\int_{a}^{b}\int_{\mathbb{R}}\phi g[u]dxdt\Bigg\}
&\leq\int_{a}^{b}\int_{\mathbb{R}}\big|g[u]\big|dxdt\\
&\leq c(b-a),
\end{align*}
where Lemma \ref{000} and Lemma \ref{00012} have been used.
\end{proof}


\section{Acknowlegement}
We would like to thank the referees who did a very careful job reading
this paper. Their indications and suggestions have 
not only contributed to make our exposition clearer, but have also
helped us to correct and improve our original results.

\end{document}